\documentclass{article}%
\usepackage{amssymb}
\usepackage{amsmath}
\usepackage{amsfonts}%
\usepackage{graphicx}
\providecommand{\U}[1]{\protect\rule{.1in}{.1in}}
\newtheorem{theorem}{Theorem}

\newtheorem{corollary}[theorem]{Corollary}

\newtheorem{definition}[theorem]{Definition}

\newtheorem{lemma}[theorem]{Lemma}

\newtheorem{proposition}[theorem]{Proposition}
\newtheorem{remark}[theorem]{Remark}

\newenvironment{proof}[1][Proof]{\noindent\textbf{#1.} }{\ \rule{0.5em}{0.5em}}
\begin{document}

\title{Saddle-point structure of fixed points in a reaction-diffusion equation with
discontinuous nonlinearity}
\author{Jos\'{e} Valero\\{\small Centro de Investigaci\'{o}n Operativa. Universidad Miguel
Hern\'{a}ndez de Elche, }\\{\small Avda. Universidad s/n, 03202-Elche (Alicante). Spain}\\{\small E.mail: jvalero@umh.es}}
\date{}
\maketitle

\begin{abstract}
In this paper, we study the local behaviour of solutions near the fixed points
of a reaction-diffusion equation with discontinuous nonlinearity. By employing
an appropriate linearization around the fixed points, which involves the Dirac
delta distribution, we analyze the stability of the stationary solutions and
demonstrate that they exhibit a saddle-point structure. As a result, we
establish the hyperbolicity of the fixed points.

\end{abstract}

\bigskip



\section{Introduction}

In this paper, we aim to study the local behaviour of solutions near the fixed
points of a reaction-diffusion equation governed by a discontinuous
nonlinearity of Heaviside type. Specifically, we consider the parabolic
equation%
\begin{equation}
u_{t}-u_{xx}\in H_{0}(u),\ x\in\left(  0,1\right)  ,\ t>0,\label{RD}%
\end{equation}
where $u$ satisfies Dirichlet boundary conditions and $H_{0}$ is the
set-valued map given in (\ref{Heaviside}). This equation arises as the limit
of a sequence of Chafee-Infante problems \cite{ArrRBVal}, a well-known and
paradigmatic model, notable for being one of the few partial differential
equations for which the structure of the global attractor has been completely
described. Furthermore, reaction-diffusion equations with this type of
nonlinearity have been used to model, for example, processes of combustion in
porous media \cite{FeNo}, the propagation of electrical impulses in nerve
axons (see \cite{Terman}, \cite{Terman2}) and the climate on Earth (see
\cite{Budyko}, \cite{NC}).

In the theory of dynamical systems, both the local and global behaviour of
solutions play essential roles. Regarding global behaviour, the existence and
structure of global attractors are key to understanding the long-term dynamics
of solutions. We emphasize that the global attractor of the Chafee-infante
problem consists of a finite number of fixed points and the heteroclinic
connections between them. Moreover, it is known exactly which connections
exist, thereby providing a complete description of the dynamics \cite{Henry85}.

It is natural to expect that the limit equation (\ref{RD}) inherits the
structure of the global attractor from the Chaffe-Infante approximations.
Problem (\ref{RD}) posseses a countable number of fixed points: $v_{0}%
=0,\ v_{1}^{\pm},\ v_{2}^{\pm},...$, which appear as a result of an infinite
cascade of bifurcations in the approximating Chafee-Infante problems.
Moreover, it is proved in \cite{ArrRBVal} that for any $n\geq1$ the sequence
of fixed points $\{v_{n,j}^{\pm}\}_{j\geq1}$ corresponding to the
approximating Chafee-Infante problems converges to $v_{n}^{\pm}$ in the space
$C^{1}([0,1])$. The conjecture about the structure of the attractor states
that there is an heteroclinic connection from a fixed point $z$ to a fixed
point $y$ (denoted $z\rightsquigarrow y$), meaning that there is a bounded
complete trajectory $\phi\left(  \text{\textperiodcentered}\right)  $ such
that $\phi\left(  t\right)  \overset{t\rightarrow+\infty}{\rightarrow}y,$
$\phi\left(  t\right)  \overset{t\rightarrow-\infty}{\rightarrow}z$, if and
only if either $z=0$, $y=v_{n}^{\pm}$, $n\geq1$, or $z=v_{n}^{\pm}, $
$y=v_{k}^{\pm}$ with $k<n$. In \cite{ArrRBVal} this question was partially
solved by showing that the above are the only possible connections and that
the following ones do occur:

\begin{enumerate}
\item $0\rightsquigarrow v_{n}^{\pm}$, $\forall n\geq1$;

\item $v_{n}^{\pm}\rightsquigarrow v_{1}^{\pm}$, if $n>2$;

\item $v_{n}^{\pm}\rightsquigarrow v_{n-1}^{\pm},\ \forall n\geq1$;\ 

\item If $v_{k}\rightsquigarrow v_{i}$, $k>i\geq1$, then $v_{nk}%
\rightsquigarrow v_{ni}$, $\forall n\geq1.$
\end{enumerate}

The existence of the other possible connections remains an open problem so far.

Another important question, which is the main focus of this paper, concerns
the local behavior of solutions near the fixed points. It was proved in
\cite{ArrRBVal} that the points $v_{1}^{\pm}$ are asymptotically stable while
the other points are unstable in the spaces $L^{2}(0,1)$ and $H_{0}^{1}(0,1)$,
which coincides with the same result for the Chafee-Infante equation. The
proof of the instability of the points $v_{n}^{\pm}$ for $n\geq2$ was simple
and based on a Lyapunov function. However, this approach does not provide
detailed information about how solutions behave when starting near the fixed
points. This is way in this paper we intend to use the classical linearization
method to both confirm the instability of the fixed points and to gain deeper
insight into the behavior of solutions to the nonlinear equation based on its
linear counterpart. In addition, we establish these results in more regular
spaces. Tipically, linearization requires smoothness of the nonlinear term in
the equation, which poses a significant challenge for equation (\ref{RD}), as
its nonlinearity is discontinuous. Inspired by the approach in \cite{Bensid},
we construct an appropriate linear approximation using the Dirac delta
distribution for any $v_{n}^{\pm}$, $n\geq2$. We begin by studying the
eigenvalue problem for the linearized equation, showing that the spectrum is
discrete, that zero is not an eigenvalue, and that the number of negative
eigenvalues is finite. In particular, for $n=2,3$ there are exactly $n$
negative eigenvalues. These results are then used to rigorously establish the
instability of the fixed points. Finally, we show that the stable and unstable
sets of the fixed points $v_{n}^{\pm}$, $n\geq2,$ are tangent to the linear
subspaces generated by the negative and positive eigenvalues of the linearized
equation, respectively. This proves that the fixed points exhibit a
saddle-point structure, and hence, are hyperbolic.

This paper is organized as follows. Section 2 introduces the problem setup. In
Section 3, we study the eigenvalue problem arising from the linearization
around the fixed points. Section 4 is devoted to the analysis of the stability
of the fixed points. Finally, in Section 5, we establish the saddle
point-property of the points $v_{n}^{\pm}$ for all $n\geq2$.

\section{Set-up of the problem}

We consider the problem%
\begin{equation}
\left\{
\begin{array}
[c]{c}%
u_{t}-u_{xx}\in H_{0}(u),\ x\in\left(  0,1\right)  ,\ t>0,\\
u(0,t)=0,\ u(1,t)=0,\\
u(x,0)=u_{0}(x),
\end{array}
\right. \label{Eq}%
\end{equation}
where%
\begin{equation}
H_{0}\left(  u\right)  =\left\{
\begin{array}
[c]{c}%
-1\text{ if }u<0,\\
\lbrack-1,1]\text{ if }u=0,\\
1\text{ if }u>0,
\end{array}
\right. \label{Heaviside}%
\end{equation}
is the Heaviside function.

It is known \cite{ArrRBVal} that (\ref{Eq}) can be rewritten in the abstract
form
\begin{equation}
\left\{
\begin{array}
[c]{l}%
\dfrac{\partial u}{\partial t}+\partial\psi^{1}(u)-\partial\psi^{2}(u)\ni0,\\
u(0)=u_{0}.
\end{array}
\right. \label{Abstract}%
\end{equation}
where $\partial\psi^{1}$ and $\partial\psi^{2}$ are the subdifferentials of
the proper, convex, lower semicontinuous funcions $\psi^{i}\colon L^{2}%
(\Omega)\rightarrow(-\infty,+\infty]$, $\Omega=\left(  0,1\right)  $, given
by
\[
\psi^{1}\left(  u\right)  =\left\{
\begin{array}
[c]{c}%
\frac{1}{2}\int_{\Omega}\left\vert \nabla u\right\vert ^{2}dx,\text{ if }u\in
H_{0}^{1}\left(  \Omega\right)  ,\\
+\infty\text{, otherwise,}%
\end{array}
\right.
\]%
\[
\partial\psi^{1}(u)=\left\{  y\in L^{2}(\Omega):y(x)=-\frac{d^{2}u}{dx^{2}%
}(x)\text{, a.e. on }\Omega\right\}  ,
\]

\[
\psi^{2}\left(  u\right)  =\left\{
\begin{array}
[c]{c}%
\int_{\Omega}\int_{0}^{u}H_{0}\left(  s\right)  dsdx,\text{ if }\int_{0}%
^{u}H_{0}\left(  s\right)  ds\in L^{1}\left(  \Omega\right)  ,\\
+\infty\text{, otherwise,}%
\end{array}
\right.
\]
and
\[
\partial\psi^{2}(u)=\left\{  y\in L^{2}\left(  \Omega\right)  :\text{
}y\left(  x\right)  \in H_{0}\left(  u\left(  x\right)  \right)  \text{, a.e.
on }\Omega\right\}  .
\]
We note that $\int_{0}^{u}H_{0}\left(  s\right)  ds=\left\vert u\right\vert $
and $D\left(  \partial\psi^{1}\right)  =H^{2}\left(  \Omega\right)  \cap
H_{0}^{1}\left(  \Omega\right)  ,\ D\left(  \partial\psi^{2}\right)
=L^{2}\left(  \Omega\right)  $. We will denote the operator $\partial\psi^{1}$
by $A.$

\begin{definition}
We say that the function $u\in C([0,\infty),L^{2}(\Omega))$ is a strong
solution of (\ref{Abstract}) if:

\begin{enumerate}
\item $u(0)=u_{0}$;

\item $u(\cdot)$ is absolutely continuous on every interval $[T_{1}%
,T_{2}]\subset(0,\infty)$ and $u\left(  t\right)  \in D\left(  \partial
\psi^{1}\right)  $ for a.a. $t$;

\item There exist a function $g\in L_{loc}^{2}(0,\infty;L^{2}\left(
\Omega\right)  )$ such that $g(t)\in\partial\psi^{2}(u(t))$, a.e. on
$(0,\infty) $, and
\begin{equation}
\frac{du(t)}{dt}+Au\left(  t\right)  -g(t)=0,\text{ for a.a. }t\in
(0,\infty).\label{exist1}%
\end{equation}

\end{enumerate}
\end{definition}

It is know (see either \cite{Valero2001} or \cite{CLV20}) that for any
$u_{0}\in L^{2}\left(  \Omega\right)  $ there exists at least one strong
solution to problem (\ref{Abstract}). In \cite{ArrRBVal} the structure of the
global attractor of the multivalued semiflow generated by the strong solutions
to problem (\ref{Abstract}) was studied. It was established that the
attractors consists of the fixed points (which are countably many) and the
heteroclinic connections between them. As described in the introduction, some
heteroclinic connections are forbidden and several ones were proved to happen.
The complete description of the connections remains an open problem so far.

In this paper, we are interested in a different task. Namely, we intend to
describe the behaviour of solutions near the fixed points, establishing that
the non-zero equilibria satisfy the saddle-point property, so they are
hyperbolic. A fixed point is a constant solution. It is known \cite{ArrRBVal}
that $v\left(  \text{\textperiodcentered}\right)  $ is a fixed point of
(\ref{Abstract}) if and only if $v\in H^{2}\left(  \Omega\right)  \cap
H_{0}^{1}\left(  \Omega\right)  $ and%
\begin{equation}
\left\{
\begin{array}
[c]{c}%
-\dfrac{d^{2}v}{dx^{2}}\in H_{0}\left(  v\right)  ,\text{ for a.a. }x\in
\Omega,\\
v\left(  0\right)  =v\left(  1\right)  =0.
\end{array}
\right. \label{Fixed}%
\end{equation}
It was proved in \cite{ArrRBVal} there is an infinite but countable number of
fixed points, which are given explicitely by:%
\begin{align}
&  v_{0}\equiv0,\label{FixedFormulas}\\
&  v_{1}^{+}(x)=-\frac{x^{2}}{2}+\frac{x}{2},\quad v_{1}^{-}(x)=-v_{1}%
(x)^{+},\nonumber\\
&  v_{2}^{+}(x)=\left\{
\begin{array}
[c]{ll}%
-\frac{x^{2}}{2}+\frac{x}{4}, & \text{if }0\leq x\leq\frac{1}{2},\\
\frac{(x-\frac{1}{2})^{2}}{2}-\frac{x-\frac{1}{2}}{4}, & \text{if }\frac{1}%
{2}\leq x\leq1,
\end{array}
\right.  \quad v_{2}^{-}(x)=-v_{2}^{+}(x),\nonumber\\
&  \vdots\nonumber\\
&  v_{n}^{+}(x)=\left\{
\begin{array}
[c]{ll}%
-\frac{x^{2}}{2}+\frac{x}{2n}, & \text{if }0\leq x\leq\frac{1}{n},\\
-\frac{(x-\frac{k}{n})^{2}}{2}+\frac{x-\frac{k}{n}}{2n}, & \text{if }\frac
{k}{n}\leq x\leq\frac{k+1}{n},k\text{ is even,}\\
\frac{(x-\frac{k}{n})^{2}}{2}-\frac{x-\frac{k}{n}}{2n} & \text{if }\frac{k}%
{n}\leq x\leq\frac{k+1}{n},k\text{ is odd,}%
\end{array}
\right.  \quad v_{n}^{-}(x)=-v_{n}^{+}(x),\ n\in\mathbb{N}.\nonumber
\end{align}

We know from \cite{ArrRBVal} that the points $v_{1}^{+},v_{1}^{-}$ are
asymptotically stable, while the other fixed points are unstable. The proof of
instability of the points $v_{n}^{+}.v_{n}^{-}.\ n\geq2,$ was given using a
Lyapunov function. However, our intention in this paper is to give a different
proof of instability, which is closer to the classical one. Namely, following
the ideas from \cite{Bensid}, we intend to prove instability using the
eigenvalues of the linearization of problem (\ref{Abstract}) on a suitable
neighborhood of the fixed points $v_{n}^{+}.v_{n}^{-}$. This is not an easy
task as we are working with a discontinuous nonlinear function.

We finish this section by proving an auxiliary result. Let $X$ be a Hilbert
space with scalar product $\left(  \text{\textperiodcentered}%
,\text{\textperiodcentered}\right)  _{X}$ and norm $\left\Vert
\text{\textperiodcentered}\right\Vert _{X}$. We consider the abstract
differential equation
\begin{equation}
\left\{
\begin{array}
[c]{c}%
\dfrac{du}{dt}+Gu=f(t),\ t>\tau,\\
u(\tau)=u_{\tau}\in X,
\end{array}
\right. \label{EqG}%
\end{equation}
where $\tau\in\mathbb{R}$, $f\in L_{loc}^{2}(\tau,\infty;X)$ and
$G:D(G)\rightarrow X$ is the infinitesimal generator of a strongly continuous
semigroup of bounded linear operators $S(t)=e^{-Gt}:X\rightarrow X$, $t\geq0
$. It is well known that $cl_{X}D(G)=X$ \cite[Corollary 2.5]{Pazy} and that
there exist $M\geq1,\ \omega\geq0$ such that $\left\Vert S(t)\right\Vert
_{\mathcal{L}(X,X)}\leq Me^{\omega t}$ for $t\geq0$ \cite[Theorem 2.2]{Pazy}.

The function $u\in C([\tau,\infty),X)$ is called a strong solution to
(\ref{EqG}) if it absolutely continuous on every compact interval of $\left(
\tau,\infty\right)  $, $u\left(  t\right)  \in D(G)$ for a.a. $t>\tau$ and the
quality in (\ref{EqG}) is satisfied for a.a. $t>\tau.$ It is called a mild
solution if%
\begin{equation}
u(t)=e^{-G(t-\tau)}u_{\tau}+\int_{\tau}^{t}e^{-G(t-s)}f(s)ds\text{, for }%
\tau\leq t.\label{Mild}%
\end{equation}
It is called a classical solution if it is continuously differentiable on
$\left(  \tau,\infty\right)  $. $u\left(  t\right)  \in D(G)$ for all $t>\tau$
and the equality in (\ref{EqG}) is satisfied for all $t>0.$

It is immediate that problem (\ref{EqG}) possesses a unique mild solution. If
$u_{\tau}\in D(A)$ and $f\left(  t\right)  $ is continuously differentiable in
$[\tau,\infty)$, then $u\left(  \text{\textperiodcentered}\right)  $ is the
unique classical solution (see \cite[p.107]{Pazy}).

\begin{lemma}
\label{MildLemma}Assume the existence of $K\geq0$ such that
\[
\left(  G(v),v\right)  _{X}\geq-K\left\Vert v\right\Vert _{X}^{2}\text{
}\forall v\in D(G).
\]
Then for any $u_{\tau}\in X$ and $f\in L_{loc}^{2}(\tau,\infty;X)$, a strong
solution $u\left(  \text{\textperiodcentered}\right)  $ to (\ref{EqG}) is also
a mild solution.
\end{lemma}

\begin{proof}
Let $u\left(  \text{\textperiodcentered}\right)  $ be a strong solution. We
take sequences $\{u_{\tau}^{n}\}\subset D\left(  G\right)  $, $\{f^{n}\left(
\text{\textperiodcentered}\right)  \}\subset C^{1}([\tau,\infty),X)$
satisfying%
\[
u_{\tau}^{n}\rightarrow u_{\tau}\text{ in }X,
\]%
\[
f^{n}\rightarrow f\text{ in }L_{loc}^{2}(\tau,\infty;X).
\]
If we replace $f,u_{\tau}$ by $f^{n},u_{\tau}^{n}$ in (\ref{EqG}), then there
is a unique classical solution $u^{n}\left(  \text{\textperiodcentered
}\right)  $.

Let $w=u^{n}-u$. By the definition of strong solution and Lemma 1.2 in
\cite[p.100]{Barbu} we obtain that%
\begin{align*}
\left\Vert w^{n}(t)\right\Vert _{X}\frac{d}{dt}\left\Vert w^{n}\right\Vert
_{X}  & =\left(  f^{n}\left(  t\right)  -f(t),w^{n}(t)\right)  _{X}-\left(
Gw^{n}(t),w^{n}(t)\right)  _{X}\\
& \leq\left\Vert f^{n}\left(  t\right)  -f(t)\right\Vert _{X}\left\Vert
w^{n}(t)\right\Vert _{X}+K\left\Vert w^{n}(t)\right\Vert _{X}^{2}\text{ for
a.a. }t>\tau.
\end{align*}
Since $\left\Vert w^{n}(t)\right\Vert _{X}$ is absolutely continuous on any
interval $[\tau+\varepsilon,T]$, $0<\varepsilon<T$, we infer using Gronwall's
lemma that%
\[
\left\Vert w^{n}(t)\right\Vert _{X}\leq e^{K(t-\tau-\varepsilon)}\left\Vert
w^{n}(\tau+\varepsilon)\right\Vert _{X}+\int_{\tau+\varepsilon}^{t}%
e^{K(t-s)}\left\Vert f^{n}\left(  s\right)  -f(s)\right\Vert _{X}ds\text{,
}t\geq\tau+\varepsilon.
\]
Passing to the limit as $\varepsilon\rightarrow0$ we have%
\[
\left\Vert w^{n}(t)\right\Vert _{X}\leq e^{K(t-\tau)}\left\Vert w^{n}%
(\tau)\right\Vert _{X}+\int_{\tau}^{t}e^{K(t-s)}\left\Vert f^{n}\left(
s\right)  -f(s)\right\Vert _{X}ds\text{, }t\geq\tau,
\]
so $u^{n}\rightarrow u$ in $C([\tau,\infty),X)$. The function $u^{n}$ is a
mild solution. Hence,%
\begin{equation}
u^{n}(t)=e^{-G(t-\tau)}u_{\tau}^{n}+\int_{\tau}^{t}e^{-G(t-s)}f^{n}%
(s)ds.\label{unMild}%
\end{equation}
By%
\[
\left\Vert \int_{\tau}^{t}e^{-G(t-s)}\left(  f^{n}(s)-f(s)\right)
ds\right\Vert _{X}\leq\int_{\tau}^{t}e^{\omega(t-s)}\left\Vert f^{n}%
(s)-f(s)\right\Vert _{X}ds\overset{n\rightarrow\infty}{\rightarrow}0
\]
and the continuity of the map $e^{-G(t-\tau)}$ we obtain passing to the limit
in (\ref{unMild}) that $u$ is a mild solution.
\end{proof}

\section{The eigenvalue problem of the linearization around fixed points}

In this section, we define first a suitable linearization of problem
(\ref{Abstract}) around the fixed points $v_{n}^{+},v_{n}^{-},\ n\geq2.$ Then
we study the eigenvalue problem for the resulting linear equation.

Let us consider first the fixed point $v_{2}^{+}$ (for $v_{2}^{-}$ the
analysis is similar). The function $v_{2}^{+}$ (see (\ref{FixedFormulas}))
belongs to $C^{1}\left(  [0,1]\right)  $. Let $u(\cdot)\in C^{1}([0,1])$ and
let $u$ vanish at $x_{1},x_{2},\ldots,x_{k}\in\left(  0,1\right)  $ and
$u^{^{\prime}}(x_{j})\neq0$. Also, let $\delta($\textperiodcentered$)$ be the
Dirac-delta function. The function $\delta\left(  \text{\textperiodcentered
}-x\right)  ,\ x\in\left(  0,1\right)  ,$ is a distribution of the space
$C_{0}^{\infty}\left(  0,1\right)  $ given by
\[
\left\langle \delta(\text{\textperiodcentered}-x),\varphi\right\rangle
=\varphi\left(  x\right)  .
\]
Then the composition $\delta(u(\cdot))$ is equal to
\[
\delta(u(\cdot))=\sum_{j=1}^{k}\frac{\delta(\cdot-x_{j})}{\left\vert
u^{^{\prime}}(x_{j})\right\vert }.
\]
We observe that if we consider the function $H_{0}\left(
\text{\textperiodcentered}\right)  $ as a distribution in $C_{0}^{\infty
}\left(  0,1\right)  $, then $H_{0}^{\prime}=2\delta\left(
\text{\textperiodcentered}\right)  $. As $\left(  v_{2}^{\pm}\right)
^{\prime}\left(  \frac{1}{2}\right)  =\mp\frac{1}{4}$ and $x=\frac{1}{2}$ is
the only point in $\left(  0,1\right)  $ where $v_{2}^{\pm}$ vanishes, the
distribution $H_{0}^{\prime}\left(  v_{2}^{\pm}\right)  $ is given by%
\[
H_{0}^{^{\prime}}\left(  v_{2}^{\pm}\right)  =2\delta\left(  v_{2}^{\pm
}\right)  =\frac{2\delta_{\frac{1}{2}}}{\left\vert \left(  v_{2}^{\pm}\right)
^{^{\prime}}\left(  \frac{1}{2}\right)  \right\vert }=8\delta_{\frac{1}{2}},
\]
where $\delta_{\frac{1}{2}}=\delta\left(  \text{\textperiodcentered}-\frac
{1}{2}\right)  .$ Thus, the linearized equation is%
\[
\dfrac{\partial U}{\partial t}-\dfrac{\partial^{2}U}{\partial x^{2}}%
=8\delta_{\frac{1}{2}}U.
\]

For any $U\in H_{0}^{1}\left(  0,1\right)  $, $8\delta_{\frac{1}{2}}U\in
H^{-1}\left(  0,1\right)  $ and%
\[
\left\langle 8\delta_{\frac{1}{2}}U,\varphi\right\rangle _{\left(
H^{-1},H_{0}^{1}\right)  }=8U\left(  \frac{1}{2}\right)  \varphi\left(
\frac{1}{2}\right)  \text{ for }\varphi\in H_{0}^{1}\left(  0,1\right)  .
\]
Thus, $H_{0}^{^{\prime}}\left(  v_{2}^{+}\right)  :H_{0}^{1}\left(
0,1\right)  \rightarrow H^{-1}\left(  0,1\right)  $ is a linear operator. We
consider the eigenvalue problem%
\begin{equation}
\left\{
\begin{array}
[c]{c}%
-U_{xx}-8\delta_{\frac{1}{2}}U=\lambda U,\\
U\left(  0\right)  =U\left(  1\right)  =0.
\end{array}
\right. \label{Eigen}%
\end{equation}
A solution of (\ref{Eigen}) is a function $U\in H_{0}^{1}\left(  0,1\right)  $
such that $U_{xx}\in H^{-1}\left(  0,1\right)  $ and the equality in
(\ref{Eigen}) is satisfied in $H^{-1}(0,1).$

\begin{proposition}
\label{EigenSol}$U\left(  \text{\textperiodcentered}\right)  \in H_{0}%
^{1}\left(  0,1\right)  $ is a solution of problem (\ref{Eigen}) if and only
if it satisfies%
\begin{equation}
\left\{
\begin{array}
[c]{c}%
-U_{xx}=\lambda U\text{ for }x\in\left(  0,\frac{1}{2}\right)  \cup\left(
\frac{1}{2},1\right)  ,\\
U\left(  0\right)  =U\left(  1\right)  =0,\ U_{-}\left(  \frac{1}{2}\right)
=U_{+}\left(  \frac{1}{2}\right)  ,\\
U_{-}^{\prime}\left(  \frac{1}{2}\right)  -U_{+}^{\prime}\left(  \frac{1}%
{2}\right)  =8U\left(  \frac{1}{2}\right)  ,
\end{array}
\right. \label{tao}%
\end{equation}
where $U_{\pm}(1/2)$ and $U_{\pm}^{^{\prime}}(1/2)$ are the directional limits
of $U$ and $U^{^{\prime}}$ on $x=1/2$.
\end{proposition}

\begin{proof}
Let $U$ be a solution to (\ref{Eigen}). The continuity property at $x=1/2$
follows from $U\in C([0,1])$. Taking an arbitrary test function $\varphi\in
H_{0}^{1}\left(  0,1\right)  $ such that $\varphi\left(  x\right)  =0$, for
$x\in\lbrack1/2,1]$ (resp. $[0,1/2]$), it is easy to see that $-U_{xx}=\lambda
U$ in $\left(  0,1/2\right)  $ (respectively, $\left(  1/2,1\right)  $).

It remains to check the jump condition on $U^{^{\prime}}$. We observe that
$U^{\prime}\in C([0,1/2])$, so $U_{-}^{\prime}\left(  1/2\right)  $ exists,
and the same is valid for $U_{+}^{\prime}\left(  1/2\right)  $. On the one
hand, we consider a test function $\zeta_{l,k}(x)\in H_{0}^{1}(0,1)$ such that%
\[
\zeta_{l,k}(x)=\left\{
\begin{array}
[c]{c}%
1\text{ if }x\in\lbrack\frac{1}{2}-l,\frac{1}{2}+l],\\
0\text{ if }x\in\lbrack0,\frac{1}{2}-l-k]\cup\lbrack\frac{1}{2}+l+k,1],\\
\zeta_{l,k}(x)\in\lbrack0,1]\text{, otherwise,}%
\end{array}
\right.
\]
where $0<k<l<\frac{1}{4}$. Also, let%
\[
\zeta_{l}(x)=\left\{
\begin{array}
[c]{c}%
1\text{ if }x\in\lbrack\frac{1}{2}-l,\frac{1}{2}+l],\\
0\text{ if }x\in\lbrack0,\frac{1}{2}-l]\cup\lbrack\frac{1}{2}+l,1].
\end{array}
\right.
\]
It is clear that $\zeta_{l,k}\rightarrow\zeta_{l}$ in $L^{2}(0,1)$ as
$k\rightarrow0$. On the other hand, we take mollifiers $h_{n}($%
\textperiodcentered$)\in C_{0}^{\infty}(\mathbb{R})$ such that $h_{n}(x)\geq
0$, $supp(h_{n})\subset\lbrack1/2-1/n,1/2+1/n]$, $\int_{1/2-1/n}%
^{1/2+1/n}h_{n}(x)dx=1$. From%
\[
\int_{0}^{1}h_{n}(x)\varphi(x)dx=\varphi(x_{n})\int_{1/2-1/n}^{1/2+1/n}%
h_{n}(x)dx\overset{n\rightarrow\infty}{\rightarrow}\varphi\left(  \frac{1}%
{2}\right)  =\left\langle \delta_{\frac{1}{2}},\varphi\right\rangle _{\left(
H^{-1},H_{0}^{1}\right)  },\text{ }\forall\varphi\in H_{0}^{1}(0,1),
\]
it follows that $h_{n}\rightarrow\delta_{\frac{1}{2}}$ weakly in
\thinspace$H^{-1}(0,1)$. We consider the approximating problems%
\begin{equation}
\left\{
\begin{array}
[c]{c}%
-(1+\alpha)U_{xx}^{n}-8h_{n}U^{n}=\lambda U-\alpha U_{xx}\text{ for }%
x\in\left(  0,1\right)  ,\\
U^{n}\left(  0\right)  =U^{n}\left(  1\right)  =0,
\end{array}
\right. \label{EigenApp}%
\end{equation}
for $\alpha>0$. We define the bilinear form $a_{n}:H_{0}^{1}(0,1)\times
H_{0}^{1}(0,1)\rightarrow\mathbb{R}$ by $a_{n}(u,v)=\left(  1+\alpha\right)
\int_{0}^{1}u_{x}v_{x}dx-8\int_{-1/2-1/n}^{1/2+1/n}h_{n}(x)u\left(  x\right)
v(x)dx$. This form is continuous and coercive for $\alpha>7$. Indeed,%
\[
\left\vert a_{n}(u,v)\right\vert \leq\left(  1+\alpha\right)  \left\Vert
u\right\Vert _{H_{0}^{1}}\left\Vert v\right\Vert _{H_{0}^{1}}+8\left\vert
u(x_{n})\right\vert \left\vert v(x_{n})\right\vert \leq C_{1}\left\Vert
u\right\Vert _{H_{0}^{1}}\left\Vert v\right\Vert _{H_{0}^{1}},
\]%
\begin{align*}
a_{n}(u,u)  & =\left(  1+\alpha\right)  \left\Vert u\right\Vert _{H_{0}^{1}%
}^{2}-8u^{2}\left(  x_{n}\right)  =\left(  1+\alpha\right)  \left\Vert
u\right\Vert _{H_{0}^{1}}^{2}-8\left(  \int_{0}^{x_{n}}u^{\prime}(x)dx\right)
^{2}\\
& \geq\left(  1+\alpha\right)  \left\Vert u\right\Vert _{H_{0}^{1}}%
^{2}-8\left\Vert u\right\Vert _{H_{0}^{1}}^{2}\left\vert x_{n}\right\vert
\geq\left(  \alpha-7\right)  \left\Vert u\right\Vert _{H_{0}^{1}}^{2}.
\end{align*}
Hence, as $g_{U}=\lambda U-\alpha U_{xx}\in H^{-1}(0,1)$, Lax-Milgram's
theorem implies the existence of a unique $U^{n}\in H_{0}^{1}\left(
0,1\right)  $ solving problem (\ref{EigenApp}). Using the last inequality we
obtain also that%
\[
\left(  \alpha-7\right)  \left\Vert U^{n}\right\Vert _{H_{0}^{1}}^{2}%
\leq\left\Vert g_{U}\right\Vert _{H^{-1}}\left\Vert U^{n}\right\Vert
_{H_{0}^{1}}\leq\frac{1}{2\left(  \alpha-7\right)  }\left\Vert g_{U}%
\right\Vert _{H^{-1}}^{2}+\frac{\alpha-7}{2}\left\Vert U^{n}\right\Vert
_{H_{0}^{1}}^{2},
\]
so $U^{n}$ is bounded in $H_{0}^{1}\left(  0,1\right)  $ and, up to a
subsequence, $U^{n}\rightarrow V$ weakly in $H_{0}^{1}\left(  0,1\right)  $,
$U_{xx}^{n}\rightarrow V_{xx}$ weakly in $H^{-1}(0,1)$. Moreover, using the
compact embedding $H_{0}^{1}(0,1)\subset C([0,1])$ we get%
\begin{align*}
\left\langle h_{n}U^{n},\varphi\right\rangle _{\left(  H^{-1},H_{0}%
^{1}\right)  }  & =\int_{-1/2-1/n}^{1/2+1/n}h_{n}(x)U^{n}\left(  x\right)
\varphi(x)dx\\
& =U^{n}(x_{n})\varphi(x_{n})\overset{n\rightarrow\infty}{\rightarrow}%
V(\frac{1}{2})\varphi(\frac{1}{2})=\left\langle \delta_{\frac{1}{2}}%
V,\varphi\right\rangle _{\left(  H^{-1},H_{0}^{1}\right)  },
\end{align*}
so $h_{n}U^{n}\rightarrow\delta_{\frac{1}{2}}V$ weakly in $H^{-1}(0,1)$.
Passing to the limit in (\ref{EigenApp}) we have%
\begin{equation}
\left\{
\begin{array}
[c]{c}%
-(1+\alpha)V_{xx}-8\delta_{\frac{1}{2}}V=\lambda U-\alpha U_{xx}\\
V(0)=V(1)=0.
\end{array}
\right. \label{EigenV}%
\end{equation}
The bilinear form $a:H_{0}^{1}(0,1)\times H_{0}^{1}(0,1)\rightarrow\mathbb{R}$
given by $a(u,v)=\left(  1+\alpha\right)  \int_{0}^{1}u_{x}v_{x}dx-8u\left(
\frac{1}{2}\right)  v(\frac{1}{2})$ is also continuous and coercive for
$\alpha>7$. Hence, problem (\ref{EigenV}) has a unique solution $V$. Since $U$
is a solution to (\ref{EigenV}), we infer that $V=U$. Hence, $U^{n}\rightarrow
U$ weakly in $H_{0}^{1}(0,1)$ and strongly in $C([0,1])$. Also, we observe
that $h_{n}$ is bounded in $L^{\infty}(0,1/2-l)$ for a fixed $l$. Hence, from
the first equality in (\ref{EigenApp}) we obtain that $U^{n}$ is bounded in
$H^{2}(0,1/2-l)$, so $U^{n}\rightarrow U$ strongly in $C^{1}([0,1/2-l])$. The
same happens in the interval $[1/2+l,1]$. Since $U^{n}\in H^{2}(0,1)$, we
obtain multiplying the equation in (\ref{EigenApp}) by $\zeta_{l,k}$ and
taking $n\geq1/l$ that%
\begin{align*}
& -(1+\alpha)\int_{\frac{1}{2}-2l}^{\frac{1}{2}+2l}U_{xx}^{n}\zeta
_{l,k}dx-8\int_{\frac{1}{2}-\frac{1}{n}}^{\frac{1}{2}+\frac{1}{n}}h_{n}%
U^{n}\zeta_{l,k}dx\\
& =\alpha\int_{\frac{1}{2}-2l}^{\frac{1}{2}+2l}U_{x}^{n}\left(  \zeta
_{l,k}\right)  _{x}dx-(U_{x}^{n}(\frac{1}{2}+l)\zeta_{l,k}(\frac{1}%
{2}+l)-U_{x}^{n}(\frac{1}{2}-l)\zeta_{l,k}(\frac{1}{2}-l))\\
& +\int_{\frac{1}{2}-l}^{\frac{1}{2}+l}U_{x}^{n}\left(  \zeta_{l,k}\right)
_{x}dx-\int_{\frac{1}{2}-l-k}^{\frac{1}{2}-l}U_{xx}^{n}\zeta_{l,k}%
dx-\int_{\frac{1}{2}+l}^{\frac{1}{2}+l+k}U_{xx}^{n}\zeta_{l,k}dx-8U^{n}%
(x_{n})\zeta_{l,k}(x_{n})\\
& =\lambda\int_{\frac{1}{2}-2l}^{\frac{1}{2}+2l}U\zeta_{l,k}dx+\alpha
\int_{\frac{1}{2}-2l}^{\frac{1}{2}+2l}U_{x}\left(  \zeta_{l,k}\right)  _{x}dx.
\end{align*}
We note that $\left(  \zeta_{l,k}\right)  _{x}(x)=0$ on $\left(
1/2-l,1/2+l\right)  $, so $\int_{\frac{1}{2}-l}^{\frac{1}{2}+l}U_{x}%
^{n}\left(  \zeta_{l,k}\right)  _{x}dx=0$. Also,
\[
\left\vert \int_{\frac{1}{2}\pm l}^{\frac{1}{2}\pm l\pm k}U_{xx}^{n}%
\zeta_{l,k}dx\right\vert \leq\left(  \int_{\frac{1}{2}\pm l}^{\frac{1}{2}\pm
l\pm k}\left(  U_{xx}^{n}\right)  ^{2}dx\right)  ^{\frac{1}{2}}k^{\frac{1}{2}%
}\rightarrow0\text{, as }k\rightarrow0.
\]
Thus, passing to the limit as $k\rightarrow0$ and using $\zeta_{l,k}(\frac
{1}{2}\pm l)=1$, $\zeta_{l,k}(x_{n})=1$ (as $x_{n}\in\lbrack
1/2-1/n,1/2+1/n]\subset\lbrack1/2-l,1/2+l]$) we obtain%
\begin{align*}
& \alpha\int_{\frac{1}{2}-2l}^{\frac{1}{2}+2l}U_{x}^{n}\left(  \zeta
_{l,k}\right)  _{x}dx+U_{x}^{n}(\frac{1}{2}-l)-U_{x}^{n}(\frac{1}{2}%
+l)-8U^{n}(x_{n})\\
& =\lambda\int_{\frac{1}{2}-2l}^{\frac{1}{2}+2l}U\zeta_{l}dx+\alpha\int%
_{\frac{1}{2}-2l}^{\frac{1}{2}+2l}U_{x}\left(  \zeta_{l,k}\right)
_{x}dx=\lambda\int_{\frac{1}{2}-l}^{\frac{1}{2}+l}Udx+\alpha\int_{\frac{1}%
{2}-2l}^{\frac{1}{2}+2l}U_{x}\left(  \zeta_{l,k}\right)  _{x}dx.
\end{align*}
Then we pass to the limit as $n\rightarrow\infty$ having that%
\[
U_{x}(\frac{1}{2}-l)-U_{x}(\frac{1}{2}+l)-8U(\frac{1}{2})=\lambda\int%
_{\frac{1}{2}-l}^{\frac{1}{2}+l}Udx.
\]
Since the last integral converges to zero as $l\rightarrow0$, we conclude the
proof of the first part.

Let $U\left(  \text{\textperiodcentered}\right)  \in H_{0}^{1}\left(
0,1\right)  $ satisfy (\ref{tao}). Then for any $\varphi\in H_{0}^{1}\left(
0,1\right)  $ we have%
\begin{align*}
\left\langle -U_{xx},\varphi\right\rangle _{\left(  H^{-1},H_{0}^{1}\right)
}  & =\int_{0}^{1}U_{x}\varphi_{x}dx=\int_{0}^{1/2}U_{x}\varphi_{x}%
dx+\int_{1/2}^{1}U_{x}\varphi_{x}dx\\
& =U_{-}^{\prime}\left(  \frac{1}{2}\right)  \varphi\left(  \frac{1}%
{2}\right)  -\int_{0}^{1/2}U_{xx}\varphi dx-U_{+}^{\prime}\left(  \frac{1}%
{2}\right)  \varphi\left(  \frac{1}{2}\right)  -\int_{1/2}^{1}U_{xx}\varphi
dx\\
& =8U(\frac{1}{2})\varphi(\frac{1}{2})+\lambda\int_{0}^{1}U\varphi dx,
\end{align*}%
\[
\left\langle \lambda U+8\delta_{\frac{1}{2}}U,\varphi\right\rangle _{\left(
H^{-1},H_{0}^{1}\right)  }=\lambda\int_{0}^{1}U\varphi dx+8(\frac{1}%
{2})\varphi(\frac{1}{2}).
\]

\end{proof}

\bigskip

The following proposition gives a complete description of the non-trivial
solutions to problem (\ref{Eigen}).

\begin{proposition}
\label{EigenvaluesU2}There is a negative eigenvalue $\lambda_{1}<0$, which is
unique. $0$ is not an eigenvalue and there is an infinite (but countable)
number of positive eigenvalues $0<\lambda_{2}\leq\lambda_{3}\leq...$
\end{proposition}

\begin{proof}
We look for a negative eigenvalue, that is, $\lambda_{1}=-\tau^{2}$,
$\tau\not =0$. According to Proposition \ref{EigenSol}, a solution $U\left(
\text{\textperiodcentered}\right)  $ to problem (\ref{Eigen}) has to satisfy
(\ref{tao}). Then%
\begin{align*}
U\left(  x\right)   & =U_{1}\left(  x\right)  =Ae^{\tau x}+Be^{-\tau x}\text{
if }x\in\left(  0,\frac{1}{2}\right)  ,\\
U\left(  x\right)   & =U_{2}\left(  x\right)  =\overline{A}e^{\tau
x}+\overline{B}e^{-\tau x}\text{ if }x\in\left(  \frac{1}{2},1\right)  .
\end{align*}
As $U\left(  0\right)  =0,\ U\left(  1\right)  =0$, we have that
\begin{align*}
A+B  & =0\Rightarrow B=-A,,\\
\overline{A}e^{\tau}+\overline{B}e^{-\tau}  & =0\Rightarrow\overline
{B}=-\overline{A}e^{2\tau}.
\end{align*}
Hence,%
\[
U_{1}\left(  x\right)  =2A\sinh\left(  \tau x\right)  ,
\]%
\[
U_{2}\left(  x\right)  =\overline{A}\left(  e^{\tau x}-e^{-\tau x+2\tau
}\right)  =\overline{A}e^{\tau}\left(  e^{-\tau\left(  1-x\right)  }%
-e^{\tau\left(  1-x\right)  }\right)  =-2\overline{A}\sinh\left(  \tau\left(
1-x\right)  \right)  .
\]
By $U_{-}\left(  \frac{1}{2}\right)  =U_{+}\left(  \frac{1}{2}\right)  $ we
get $A=-\overline{A}$. Finally, $U_{-}^{\prime}\left(  \frac{1}{2}\right)
-U_{+}^{\prime}\left(  \frac{1}{2}\right)  =8U\left(  \frac{1}{2}\right)  $
implies that $4A\tau\cosh\left(  \frac{\tau}{2}\right)  =16A\sinh\left(
\frac{\tau}{2}\right)  .$ The equation%
\[
\tau=4\tanh\left(  \frac{\tau}{2}\right)
\]
posseses two solutions $\tau=\tau^{\ast}>0$, $\tau=-\tau^{\ast}$, which give
rise to a unique \ negative eigenvalue $\lambda_{1}=-\left(  \tau^{\ast
}\right)  ^{2}.$

Let us check that $\lambda=0$ is not an eigenvalue. From (\ref{tao}) with
$\tau=0$ we have that%
\begin{align*}
U\left(  x\right)   & =U_{1}\left(  x\right)  =A+Bx\text{ if }x\in\left(
0,\frac{1}{2}\right)  ,\\
U\left(  x\right)   & =U_{2}\left(  x\right)  =\overline{A}+\overline
{B}x\text{ if }x\in\left(  \frac{1}{2},1\right)  .
\end{align*}
Since $U\left(  0\right)  =U\left(  1\right)  =0$, we have that
$A=0,\ \overline{B}=-\overline{A}$. Further, $U_{-}\left(  \frac{1}{2}\right)
=U_{+}\left(  \frac{1}{2}\right)  $ implies that $B=\overline{A}$. Finally,
$U_{-}^{\prime}\left(  \frac{1}{2}\right)  -U_{+}^{\prime}\left(  \frac{1}%
{2}\right)  =8U\left(  \frac{1}{2}\right)  $ gives $2B=4B$, so $B=0$. Hence,
$U\left(  x\right)  \equiv0.$

Finally, if $\lambda=\tau^{2}$, $\tau\not =0$, then%
\begin{align*}
U\left(  x\right)   & =U_{1}\left(  x\right)  =A\cos\left(  \tau x\right)
+B\sin\left(  \tau x\right)  \text{ if }x\in\left(  0,\frac{1}{2}\right)  ,\\
U\left(  x\right)   & =U_{2}\left(  x\right)  =\overline{A}\cos\left(  \tau
x\right)  +\overline{B}\sin\left(  \tau x\right)  \text{ if }x\in\left(
\frac{1}{2},1\right)  .
\end{align*}
From $U\left(  0\right)  =U\left(  1\right)  =0$, we have $A=0$ and
$\overline{A}=-\overline{B}\tan(\tau)$. Hence,%
\[
U_{1}\left(  x\right)  =B\sin\left(  \tau x\right)  ,
\]%
\[
U_{2}\left(  x\right)  =\frac{\overline{B}}{\cos\left(  \tau\right)  }\left(
\sin\left(  \tau x\right)  \cos\left(  \tau\right)  -\cos\left(  \tau
x\right)  \sin\left(  \tau\right)  \right)  =C\sin\left(  \tau\left(
1-x\right)  \right)  .
\]
By $U_{-}\left(  \frac{1}{2}\right)  =U_{+}\left(  \frac{1}{2}\right)  $ we
get $B=C$. Finally, $U_{-}^{\prime}\left(  \frac{1}{2}\right)  -U_{+}^{\prime
}\left(  \frac{1}{2}\right)  =8U\left(  \frac{1}{2}\right)  $ gives
$2B\tau\cos\left(  \frac{\tau}{2}\right)  =8B\sin\left(  \frac{\tau}%
{2}\right)  $. The equation%
\[
\tau=4\tan\left(  \frac{\tau}{2}\right)
\]
has in the intervals $(\left(  2k-3\right)  \pi,\left(  2k-1\right)
\pi),\ (-\left(  2k-1\right)  \pi,-\left(  2k-3\right)  \pi),\ k=2,3,...$, two
solutions $\tau=\tau_{k}>0,\ \tau=-\tau_{k}$, which give rise to the sequence
of positive eigenvalues $\lambda_{k}=\tau_{k}^{2}$, $k=2,3,...$
\end{proof}

\bigskip

Let us consider now the fixed points $v_{n}^{\pm},\ n\geq2$. As $\left(
v_{n}^{\pm}\right)  ^{\prime}\left(  \frac{k}{n}\right)  =\pm\frac{\left(
-1\right)  ^{k}}{2n}$, $1\leq k\leq n-1,$ and $x_{k}=\frac{k}{n},$ $1\leq
k\leq n-1,$ are the points in $\left(  0,1\right)  $ where $v_{n}^{+}$
vanishes, the distribution $H_{0}^{\prime}\left(  v_{n}^{\pm}\right)  $ is
given by%
\[
H_{0}^{^{\prime}}\left(  v_{n}^{\pm}\right)  =2\delta\left(  v_{n}^{\pm
}\right)  =2\sum_{k=1}^{n-1}\frac{\delta_{x_{k}}}{\left\vert \left(
v_{n}^{\pm}\right)  ^{\prime}\left(  x_{k}\right)  \right\vert }=4n\sum
_{k=1}^{n-1}\delta_{x_{k}},
\]
where $\delta_{x_{k}}=\delta\left(  \text{\textperiodcentered}-x_{k}\right)
.$ Thus, the linearized equation is%
\[
\dfrac{\partial U}{\partial t}-\dfrac{\partial^{2}U}{\partial x^{2}}%
=4n\sum_{k=1}^{n-1}\delta_{x_{k}}U.
\]

For any $U\in H_{0}^{1}\left(  0,1\right)  $, $4n\sum_{k=1}^{n-1}\left(
-1\right)  ^{k}\delta_{x_{k}}U\in H^{-1}\left(  0,1\right)  $ and%
\[
\left\langle 4n\sum_{k=1}^{n-1}\delta_{x_{k}}U,\varphi\right\rangle _{\left(
H^{-1},H_{0}^{1}\right)  }=4n\sum_{k=1}^{n-1}U\left(  x_{k}\right)
\varphi\left(  x_{k}\right)  \text{ for }\varphi\in H_{0}^{1}\left(
0,1\right)  .
\]
Thus, $H_{0}^{^{\prime}}\left(  v_{n}^{\pm}\right)  :H_{0}^{1}\left(
0,1\right)  \rightarrow H^{-1}\left(  0,1\right)  $ is a linear operator. As
before, we consider the eigenvalue problem%
\begin{equation}
\left\{
\begin{array}
[c]{c}%
-U_{xx}-4n\sum_{k=1}^{n-1}\delta_{x_{k}}U=\lambda U,\\
U\left(  0\right)  =U\left(  1\right)  =0.
\end{array}
\right. \label{Eigen2}%
\end{equation}
A solution of (\ref{Eigen2}) is a function $U\in H_{0}^{1}\left(  0,1\right)
$ such that $U_{xx}\in H^{-1}\left(  0,1\right)  $ and the equality in
(\ref{Eigen2}) is satisfied in $H^{-1}(0,1).$

\begin{proposition}
\label{EigenSoln}$U\left(  \text{\textperiodcentered}\right)  \in H_{0}%
^{1}\left(  0,1\right)  $ is a solution of problem (\ref{Eigen2}) if and only
if it satisfies%
\begin{equation}
\left\{
\begin{array}
[c]{c}%
-U_{xx}=\lambda U\text{ for }x\in\left(  \frac{k}{n},\frac{k+1}{n}\right)
,\ k=0,...,n-1,\\
U\left(  0\right)  =U\left(  1\right)  =0,\ \\
U_{-}\left(  \frac{k}{n}\right)  =U_{+}\left(  \frac{k}{n}\right)
,\ \ k=1,...,n-1,\\
U_{-}^{\prime}\left(  \frac{k}{n}\right)  -U_{+}^{\prime}\left(  \frac{k}%
{n}\right)  =4nU\left(  \frac{k}{n}\right)  ,\ \ k=1,...,n-1,
\end{array}
\right. \label{taon}%
\end{equation}
where $U_{\pm}(\frac{k}{n})$ and $U_{\pm}^{^{\prime}}(\frac{k}{n})$ are the
directional limits of $U$ and $U^{^{\prime}}=U_{x}$ on $x_{k}=\frac{k}{n}$.
\end{proposition}

\begin{proof}
Let $U$ be a solution to (\ref{Eigen2}). The continuity property at
$x=\frac{k}{n}$ follows from $U\in C([0,1])$. Taking an arbitrary test
function $\varphi\in H_{0}^{1}\left(  0,1\right)  $ such that $\varphi\left(
x\right)  =0$, $x\in\lbrack0,1]\backslash(\frac{k}{n},\frac{k+1}{n})$, it is
easy to see that $-U_{xx}=\lambda U$ in $\left(  \frac{k}{n},\frac{k+1}%
{n}\right)  $. It remains to check the jump condition on $U^{^{\prime}}$. We
observe that $U^{\prime}\in C([\frac{k}{n},\frac{k+1}{n}])$, for any
$k=0,...,n-1$, so $U_{-}^{\prime}\left(  \frac{k}{n}\right)  ,U_{+}^{\prime
}\left(  \frac{k}{n}\right)  $ exist for $k=1,...,n-1$. On the one hand, fix
$\overline{k}\in\{1,...,n-1\}$ and consider a test function $\zeta_{l,j}\in
H_{0}^{1}(0,1)$ such that%
\[
\zeta_{l,j}(x)=\left\{
\begin{array}
[c]{c}%
1\text{ if }x\in\lbrack\frac{\overline{k}}{n}-l,\frac{\overline{k}}{n}+l],\\
0\text{ if }x\in\lbrack0,\frac{\overline{k}}{n}-l-j]\cup\lbrack\frac
{\overline{k}}{n}+l+j,1],\\
\zeta_{l,j}(x)\in\lbrack0,1]\text{, otherwise,}%
\end{array}
\right.
\]
where $0<j<l<\frac{1}{4n}$. Also, let%
\[
\zeta_{l}(x)=\left\{
\begin{array}
[c]{c}%
1\text{ if }x\in\lbrack\frac{\overline{k}}{n}-l,\frac{\overline{k}}{n}+l],\\
0\text{ if }x\in\lbrack0,\frac{\overline{k}}{n}-l]\cup\lbrack\frac
{\overline{k}}{n}+l,1].
\end{array}
\right.
\]
It is clear that $\zeta_{l,j}\rightarrow\zeta_{l}$ in $L^{2}(0,1)$ as
$j\rightarrow0$. On the other hand, for $k\in\{1,...,n-1\}$ we take mollifiers
$h_{m,k}($\textperiodcentered$)\in C_{0}^{\infty}(\mathbb{R})$ such that
$h_{m,k}(x)\geq0$, $supp(h_{m,k})\subset\lbrack\frac{k}{n}-\frac{1}{m}%
,\frac{k}{n}+\frac{1}{m}]$, $\int_{k/n-1/m}^{k/n+1/m}h_{m}(x)dx=1$. From%
\[
\int_{0}^{1}h_{m,k}(x)\varphi(x)dx=\varphi(x_{m})\int_{k/n-1/m}^{k/n+1/m}%
h_{m,k}(x)dx\overset{m\rightarrow\infty}{\rightarrow}\varphi\left(  \frac
{k}{n}\right)  =\left\langle \delta_{\frac{k}{n}},\varphi\right\rangle
_{\left(  H^{-1},H_{0}^{1}\right)  },\text{ }\forall\varphi\in H_{0}^{1}(0,1),
\]
it follows that $h_{m,k}\underset{m\rightarrow\infty}{\rightarrow}%
\delta_{\frac{k}{n}}$ weakly in \thinspace$H^{-1}(0,1)$. We consider the
approximating problems
\begin{equation}
\left\{
\begin{array}
[c]{c}%
-(1+\alpha)U_{xx}^{m}-4n\sum_{k=1}^{n-1}h_{m,k}U^{m}=\lambda U-\alpha
U_{xx}\text{ for }x\in\left(  0,1\right)  ,\\
U^{m}\left(  0\right)  =U^{m}\left(  1\right)  =0,
\end{array}
\right. \label{EigenAppn}%
\end{equation}
for $\alpha>0$. We define the bilinear form $a_{m}:H_{0}^{1}(0,1)\times
H_{0}^{1}(0,1)\rightarrow\mathbb{R}$ by $a_{m}(u,v)=\left(  1+\alpha\right)
\int_{0}^{1}u_{x}v_{x}dx-4n\sum_{k=1}^{n-1}\int_{k/n-1/m}^{k/n+1/m}%
h_{m,k}(x)u\left(  x\right)  v(x)dx$. This form is continuous and coercive for
$\alpha>4n(n-1)-1$. Indeed,%
\[
\left\vert a_{m}(u,v)\right\vert \leq\left(  1+\alpha\right)  \left\Vert
u\right\Vert _{H_{0}^{1}}\left\Vert v\right\Vert _{H_{0}^{1}}+4n\sum
_{k=1}^{n-1}\left\vert u(x_{m,k})\right\vert \left\vert v(x_{m,k})\right\vert
\leq C_{1}\left\Vert u\right\Vert _{H_{0}^{1}}\left\Vert v\right\Vert
_{H_{0}^{1}},
\]%
\begin{align*}
a_{m}(u,u)  & =\left(  1+\alpha\right)  \left\Vert u\right\Vert _{H_{0}^{1}%
}^{2}-4n\sum_{k=1}^{n-1}u^{2}\left(  x_{m,k}\right)  =\left(  1+\alpha\right)
\left\Vert u\right\Vert _{H_{0}^{1}}^{2}-4n\sum_{k=1}^{n-1}\left(  \int%
_{0}^{x_{m,k}}u^{\prime}(x)dx\right)  ^{2}\\
& \geq\left(  1+\alpha\right)  \left\Vert u\right\Vert _{H_{0}^{1}}%
^{2}-4n\left\Vert u\right\Vert _{H_{0}^{1}}^{2}\sum_{k=1}^{n-1}\left\vert
x_{m,k}\right\vert \geq\left(  \alpha-4n(n-1)+1\right)  \left\Vert
u\right\Vert _{H_{0}^{1}}^{2}.
\end{align*}
Hence, as $g_{U}=\lambda U-\alpha U_{xx}\in H^{-1}(0,1)$, Lax-Milgram's
theorem implies the existence of a unique $U^{m}\in H_{0}^{1}\left(
0,1\right)  $ solving problem (\ref{EigenAppn}). Using the last inequality we
obtain also that%
\[
\left(  \alpha-4n(n-1)+1\right)  \left\Vert U^{m}\right\Vert _{H_{0}^{1}}%
^{2}\leq\left\Vert g_{U}\right\Vert _{H^{-1}}\left\Vert U^{m}\right\Vert
_{H_{0}^{1}}\leq\frac{1}{2\left(  \alpha-4n(n-1)+1\right)  }\left\Vert
g_{U}\right\Vert _{H^{-1}}^{2}+\frac{\alpha-4n(n-1)+1}{2}\left\Vert
U^{n}\right\Vert _{H_{0}^{1}}^{2},
\]
so $U^{m}$ is bounded in $H_{0}^{1}\left(  0,1\right)  $ and, up to a
subsequence, $U^{m}\rightarrow V$ weakly in $H_{0}^{1}\left(  0,1\right)  $,
$U_{xx}^{m}\rightarrow V_{xx}$ weakly in $H^{-1}(0,1)$. Moreover, using the
compact embedding $H_{0}^{1}(0,1)\subset C([0,1])$ we get%
\begin{align*}
\left\langle \sum_{k=1}^{n-1}h_{m,k}U^{m},\varphi\right\rangle _{\left(
H^{-1},H_{0}^{1}\right)  }  & =\sum_{k=1}^{n-1}\int_{k/n-1/m}^{k/n+1/m}%
h_{m,k}(x)U^{m}\left(  x\right)  \varphi(x)dx\\
& =\sum_{k=1}^{n-1}U^{m}(x_{m,k})\varphi(x_{m,k})\overset{m\rightarrow
\infty}{\rightarrow}\sum_{k=1}^{n-1}V(\frac{k}{n})\varphi(\frac{k}{n}%
)=\sum_{k=1}^{n-1}\left\langle \delta_{\frac{k}{n}}V,\varphi\right\rangle
_{\left(  H^{-1},H_{0}^{1}\right)  },
\end{align*}
so $\sum_{k=1}^{n-1}h_{m,k}U^{m}\rightarrow\sum_{k=1}^{n-1}\delta_{\frac{k}%
{n}}V$ weakly in $H^{-1}(0,1)$. Passing to the limit in (\ref{EigenAppn}) we
have%
\begin{equation}
\left\{
\begin{array}
[c]{c}%
-(1+\alpha)V_{xx}-4n\sum_{k=1}^{n-1}\delta_{\frac{k}{n}}V=\lambda U-\alpha
U_{xx}\\
V(0)=V(1)=0.
\end{array}
\right. \label{EigenVn}%
\end{equation}
The bilinear form $a:H_{0}^{1}(0,1)\times H_{0}^{1}(0,1)\rightarrow\mathbb{R}$
given by $a(u,v)=\left(  1+\alpha\right)  \int_{0}^{1}u_{x}v_{x}%
dx-4n\sum_{k=1}^{n-1}u\left(  \frac{k}{n}\right)  v(\frac{k}{n})$ is also
continuous and coercive for $\alpha>4n(n-1)-1$. Hence, problem (\ref{EigenVn})
has a unique solution $V$. Since $U$ is a solution to (\ref{EigenVn}), we
infer that $V=U$. Hence, $U^{m}\rightarrow U$ weakly in $H_{0}^{1}(0,1)$ and
strongly in $C([0,1])$. Also, we observe that $\sum_{k=1}^{n-1}h_{m,k}$ is
bounded in $L^{\infty}(\frac{k-1}{n}+l,\frac{k}{n}-l)$ for a fixed $l$ and
$k\in\{1,...,n\}$. Hence, from the first equality in (\ref{EigenAppn}) we
obtain that $U^{m}$ is bounded in $H^{2}(\frac{k-1}{n}+l,\frac{k}{n}-l)$, so
$U^{m}\rightarrow U$ strongly in $C^{1}([\frac{k-1}{n}+l,\frac{k}{n}-l])$.
Since $U^{m}\in H^{2}(0,1)$, we obtain multiplying the equation in
(\ref{EigenAppn}) by $\zeta_{l,j}$ and taking $m\geq1/l$ that%
\begin{align*}
& -(1+\alpha)\int_{\frac{\overline{k}}{n}-2l}^{\frac{\overline{k}}{n}%
+2l}U_{xx}^{m}\zeta_{l,j}dx-4n\int_{\frac{\overline{k}}{n}-\frac{1}{m}}%
^{\frac{\overline{k}}{n}+\frac{1}{m}}h_{m,\overline{k}}U^{m}\zeta_{l,j}dx\\
& =\alpha\int_{\frac{\overline{k}}{n}-2l}^{\frac{\overline{k}}{n}+2l}U_{x}%
^{m}\left(  \zeta_{l,j}\right)  _{x}dx-(U_{x}^{m}(\frac{\overline{k}}%
{n}+l)\zeta_{l,j}(\frac{\overline{k}}{n}+l)-U_{x}^{m}(\frac{\overline{k}}%
{n}-l)\zeta_{l,j}(\frac{\overline{k}}{n}-l))\\
& +\int_{\frac{\overline{k}}{n}-l}^{\frac{\overline{k}}{n}+l}U_{x}^{m}\left(
\zeta_{l,j}\right)  _{x}dx-\int_{\frac{\overline{k}}{n}-l-j}^{\frac
{\overline{k}}{n}-l}U_{xx}^{m}\zeta_{l,j}dx-\int_{\frac{\overline{k}}{n}%
+l}^{\frac{\overline{k}}{n}+l+j}U_{xx}^{m}\zeta_{l,j}dx-4nU^{m}(x_{m}%
)\zeta_{l,j}(x_{m})\\
& =\lambda\int_{\frac{\overline{k}}{n}-2l}^{\frac{\overline{k}}{n}+2l}%
U\zeta_{l,j}dx+\alpha\int_{\frac{\overline{k}}{n}-2l}^{\frac{\overline{k}}%
{n}+2l}U_{x}\left(  \zeta_{l,j}\right)  _{x}dx.
\end{align*}
We note that $\left(  \zeta_{l,j}\right)  _{x}(x)=0$ on $\left(
\frac{\overline{k}}{n}-l,\frac{\overline{k}}{n}+l\right)  $, so $\int%
_{\frac{\overline{k}}{n}-l}^{\frac{\overline{k}}{n}+l}U_{x}^{n}\left(
\zeta_{l,j}\right)  _{x}dx=0$. Also,
\[
\left\vert \int_{\frac{\overline{k}}{n}\pm l}^{\frac{\overline{k}}{n}\pm l\pm
j}U_{xx}^{m}\zeta_{l,j}dx\right\vert \leq\left(  \int_{\frac{\overline{k}}%
{n}\pm l}^{\frac{\overline{k}}{n}\pm l\pm j}\left(  U_{xx}^{m}\right)
^{2}dx\right)  ^{\frac{1}{2}}j^{\frac{1}{2}}\rightarrow0\text{, as
}j\rightarrow0.
\]
Thus, passing to the limit as $j\rightarrow0$ and using $\zeta_{l,k}%
(\frac{\overline{k}}{n}\pm l)=1$, $\zeta_{l,k}(x_{m})=1$ (as $x_{m}\in
\lbrack\frac{\overline{k}}{n}-\frac{1}{m},\frac{\overline{k}}{n}+\frac{1}%
{m}]\subset\lbrack\frac{\overline{k}}{n}-l,\frac{\overline{k}}{n}+l]$) we
obtain%
\begin{align*}
& \alpha\int_{\frac{\overline{k}}{n}-2l}^{\frac{\overline{k}}{n}+2l}U_{x}%
^{m}\left(  \zeta_{l,j}\right)  _{x}dx+U_{x}^{m}(\frac{\overline{k}}%
{n}-l)-U_{x}^{m}(\frac{\overline{k}}{n}+l)-4nU^{m}(x_{m})\\
& =\lambda\int_{\frac{\overline{k}}{n}-2l}^{\frac{\overline{k}}{n}+2l}%
U\zeta_{l}dx+\alpha\int_{\frac{\overline{k}}{n}-2l}^{\frac{\overline{k}}%
{n}+2l}U_{x}\left(  \zeta_{l,j}\right)  _{x}dx=\lambda\int_{\frac{\overline
{k}}{n}-l}^{\frac{\overline{k}}{n}+l}Udx+\alpha\int_{\frac{\overline{k}}%
{n}-2l}^{\frac{\overline{k}}{n}+2l}U_{x}\left(  \zeta_{l,j}\right)  _{x}dx.
\end{align*}
Then we pass to the limit as $m\rightarrow\infty$ having that%
\[
U_{x}(\frac{\overline{k}}{n}-l)-U_{x}(\frac{\overline{k}}{n}+l)-4nU(\frac
{\overline{k}}{n})=\lambda\int_{\frac{1}{2}-l}^{\frac{1}{2}+l}Udx.
\]
Since the last integral converges to zero as $l\rightarrow0$, we conclude the
proof of the first part.

Let $U\left(  \text{\textperiodcentered}\right)  \in H_{0}^{1}\left(
0,1\right)  $ satisfy (\ref{taon}). Then for any $\varphi\in H_{0}^{1}\left(
0,1\right)  $ we have%
\begin{align*}
\left\langle -U_{xx},\varphi\right\rangle _{\left(  H^{-1},H_{0}^{1}\right)
}  & =\int_{0}^{1}U_{x}\varphi_{x}dx=\sum_{k=1}^{n}\int_{\frac{k-1}{n}}%
^{\frac{k}{n}}U_{x}\varphi_{x}dx\\
& =\sum_{k=1}^{n-1}\left(  U_{-}^{\prime}(\frac{k}{n})-U_{+}^{\prime}(\frac
{k}{n})\right)  \varphi(\frac{k}{n})-\sum_{k=1}^{n}\int_{\frac{k-1}{n}}%
^{\frac{k}{n}}U_{xx}\varphi dx\\
& =4n\sum_{k=1}^{n-1}U(\frac{k}{n})\varphi(\frac{k}{n})+\lambda\int_{0}%
^{1}U\varphi dx,
\end{align*}%
\[
\left\langle \lambda U+4n\sum_{k=1}^{n-1}\delta_{x_{k}}U,\varphi\right\rangle
_{\left(  H^{-1},H_{0}^{1}\right)  }=\lambda\int_{0}^{1}U\varphi
dx+4n\sum_{k=1}^{n-1}U(\frac{k}{n})\varphi(\frac{k}{n}).
\]

\end{proof}

\bigskip

Let us consider now the eigenvalue problem (\ref{Eigen2}) for $n=3:$%
\begin{equation}
\left\{
\begin{array}
[c]{c}%
-U_{xx}-12(\delta_{\frac{1}{3}}+\delta_{\frac{2}{3}})U=\lambda U,\\
U\left(  0\right)  =U\left(  1\right)  =0,
\end{array}
\right. \label{Eigen3}%
\end{equation}
which corresponds to the linearization around the fixed point $v_{3}^{+}$.

The following proposition gives a complete description of the non-trivial
solutions to problem (\ref{Eigen3}).

\begin{proposition}
\label{EigenvaluesU3}There are exactly two negative eigenvalues $\lambda
_{1}<\lambda_{2}<0$. $0$ is not an eigenvalue and there is an infinite (but
countable) number of positive eigenvalues $0<\lambda_{3}\leq\lambda_{4}... $
\end{proposition}

\begin{proof}
First, we look for negative eigenvalues $\lambda=-\tau^{2}$, $\tau\not =0$.
According to Proposition \ref{EigenSoln}, a solution $U\left(
\text{\textperiodcentered}\right)  $ to problem (\ref{Eigen3}) has to satisfy
\begin{align*}
U\left(  x\right)   & =U_{1}\left(  x\right)  =Ae^{\tau x}+Be^{-\tau x}\text{
if }x\in\left(  0,\frac{1}{3}\right)  ,\\
U\left(  x\right)   & =U_{2}\left(  x\right)  =Ce^{\tau x}+De^{-\tau x}\text{
if }x\in\left(  \frac{1}{3},\frac{2}{3}\right)  ,\\
U\left(  x\right)   & =U_{3}\left(  x\right)  =Ee^{\tau x}+Fe^{-\tau x}\text{
if }x\in\left(  \frac{2}{3},1\right)  .
\end{align*}
Since $U\left(  0\right)  =0,\ U\left(  1\right)  =0$, we obtain
\begin{align*}
A+B  & =0\Rightarrow B=-A,\\
Ee^{\tau}+Fe^{-\tau}  & =0\Rightarrow F=-Ee^{2\tau}.
\end{align*}
Hence,%
\[
U_{1}(x)=2A\sinh(\tau x),
\]%
\[
U_{3}\left(  x\right)  =E\left(  e^{\tau x}-e^{-\tau x+2\tau}\right)
=Ee^{\tau}\left(  e^{-\tau\left(  1-x\right)  }-e^{\tau\left(  1-x\right)
}\right)  =-2E\sinh\left(  \tau\left(  1-x\right)  \right)  .
\]
The continuity conditions imply that%
\[
U_{1}(\frac{1}{3})=2A\sinh(\frac{\tau}{3})=U_{2}\left(  \frac{1}{3}\right)
=Ce^{\frac{\tau}{3}}+De^{-\frac{\tau}{3}},
\]%
\[
A=\frac{Ce^{\frac{\tau}{3}}+De^{-\frac{\tau}{3}}}{2\sinh(\frac{\tau}{3})},
\]%
\[
U_{2}\left(  \frac{2}{3}\right)  =Ce^{\frac{2\tau}{3}}+De^{-\frac{2\tau}{3}%
}=U_{3}\left(  \frac{2}{3}\right)  =-2E\sinh\left(  \tau\left(  1-\frac{2}%
{3}\right)  \right)  ,
\]%
\[
E=-\frac{Ce^{\frac{2\tau}{3}}+De^{-\frac{2\tau}{3}}}{2\sinh\left(  \frac{\tau
}{3}\right)  }.
\]
Thus,%
\begin{align*}
U_{1}(x)  & =\frac{Ce^{\frac{\tau}{3}}+De^{-\frac{\tau}{3}}}{\sinh(\frac{\tau
}{3})}\sinh(\tau x),\\
U_{3}(x)  & =\frac{Ce^{\frac{2\tau}{3}}+De^{-\frac{2\tau}{3}}}{\sinh\left(
\frac{\tau}{3}\right)  }\sinh\left(  \tau\left(  1-x\right)  \right)  .
\end{align*}
By the jump conditions we have%
\begin{align*}
U_{1}^{\prime}\left(  \frac{1}{3}\right)  -U_{2}^{\prime}\left(  \frac{1}%
{3}\right)   & =\tau\left(  \frac{Ce^{\frac{\tau}{3}}+De^{-\frac{\tau}{3}}%
}{\sinh(\frac{\tau}{3})}\cosh(\frac{\tau}{3})-Ce^{\frac{\tau}{3}}%
+De^{-\frac{\tau}{3}}\right) \\
& =12U_{1}\left(  \frac{1}{3}\right)  =12(Ce^{\frac{\tau}{3}}+De^{-\frac{\tau
}{3}}),
\end{align*}%
\begin{align*}
U_{2}^{\prime}\left(  \frac{2}{3}\right)  -U_{3}^{\prime}\left(  \frac{2}%
{3}\right)   & =\tau\left(  Ce^{\frac{2\tau}{3}}-De^{-\frac{2\tau}{3}}%
+\frac{Ce^{\frac{2\tau}{3}}+De^{-\frac{2\tau}{3}}}{\sinh\left(  \frac{\tau}%
{3}\right)  }\cosh\left(  \frac{\tau}{3}\right)  \right) \\
& =12U_{2}\left(  \frac{2}{3}\right)  =12\left(  Ce^{\frac{2\tau}{3}%
}+De^{-\frac{2\tau}{3}}\right)  .
\end{align*}
Then%
\[
Ce^{\frac{\tau}{3}}\left(  \tau\frac{1}{\tanh\left(  \frac{\tau}{3}\right)
}-\tau-12\right)  =-De^{-\frac{\tau}{3}}\left(  \tau\frac{1}{\tanh\left(
\frac{\tau}{3}\right)  }+\tau-12\right)  ,
\]%
\[
Ce^{\frac{2\tau}{3}}\left(  \tau\frac{1}{\tanh\left(  \frac{\tau}{3}\right)
}+\tau-12\right)  =-De^{-\frac{2\tau}{3}}\left(  \tau\frac{1}{\tanh\left(
\frac{\tau}{3}\right)  }-\tau-12\right)  ,
\]
so%
\[
\tau\frac{1}{\tanh\left(  \frac{\tau}{3}\right)  }-\tau-12=\pm e^{\frac{\tau
}{3}}\left(  \tau\frac{1}{\tanh\left(  \frac{\tau}{3}\right)  }+\tau
-12\right)  .
\]
The solutions of the equation are the zeros of the maps
\begin{align*}
f\left(  \tau\right)   & =\tau-\left(  \tau+12\right)  \tanh\left(  \frac
{\tau}{3}\right)  -e^{\frac{\tau}{3}}\left(  \tau+\left(  \tau-12\right)
\tanh\left(  \frac{\tau}{3}\right)  \right)  ,\\
g\left(  \tau\right)   & =\tau-\left(  \tau+12\right)  \tanh\left(  \frac
{\tau}{3}\right)  +e^{\frac{\tau}{3}}\left(  \tau+\left(  \tau-12\right)
\tanh\left(  \frac{\tau}{3}\right)  \right)  ,
\end{align*}
each of which has exactly one positive root, denoted by $\tau_{f}$ and
$\tau_{g}$, respectively. Moreover, $\tau_{f}<\tau_{g}$ and, since
$f(-\tau)=e^{-\frac{\tau}{3}}f\left(  \tau\right)  $, the negative roots are
$-\tau_{f}$ and $-\tau_{g}$. Thus, there are two negative eigenvales
$\lambda_{1}=-\tau_{g}^{2}<\lambda_{2}=-\tau_{f}^{2}$.

We will prove further that $\lambda=0$ is not an eigenvalue. From (\ref{tao})
with $\tau=0$ we have that%
\begin{align*}
U\left(  x\right)   & =U_{1}\left(  x\right)  =A+Bx\text{ if }x\in\left(
0,\frac{1}{3}\right)  ,\\
U\left(  x\right)   & =U_{2}\left(  x\right)  =C+Dx\text{ if }x\in\left(
\frac{1}{3},\frac{2}{3}\right)  ,\\
U\left(  x\right)   & =U_{3}\left(  x\right)  =E+Fx\text{ if }x\in\left(
\frac{2}{3},1\right)  .
\end{align*}
By $U\left(  0\right)  =U\left(  1\right)  =0$ we have $A=0,F=-E$. The
continuity at $\frac{1}{3}$ and $\frac{2}{3}$ gives $B=3C+D,\ E=3C+2D.$
Finally,%
\[
U_{1}^{\prime}\left(  \frac{1}{3}\right)  -U_{2}^{\prime}\left(  \frac{1}%
{3}\right)  =3C=12U_{1}\left(  \frac{1}{3}\right)  =12C+4D,
\]%
\[
U_{2}^{\prime}\left(  \frac{2}{3}\right)  -U_{3}^{\prime}\left(  \frac{2}%
{3}\right)  =3C+3D=12U_{2}\left(  \frac{2}{3}\right)  =12C+8D,
\]
so $C=D=B=F=0.$

If $\lambda=\tau^{2}$, $\tau\not =0$, then%
\begin{align*}
U\left(  x\right)   & =U_{1}\left(  x\right)  =A\cos\left(  \tau x\right)
+B\sin\left(  \tau x\right)  \text{ if }x\in\left(  0,\frac{1}{3}\right)  ,\\
U(x)  & =U_{2}\left(  x\right)  =C\cos\left(  \tau x\right)  +D\sin\left(
\tau x\right)  \text{ if }x\in\left(  \frac{1}{3},\frac{2}{3}\right)  ,\\
U\left(  x\right)   & =U_{3}\left(  x\right)  =E\cos\left(  \tau x\right)
+F\sin\left(  \tau x\right)  \text{ if }x\in\left(  \frac{2}{3},1\right)  .
\end{align*}
From $U_{1}\left(  0\right)  =0$, $U_{3}\left(  1\right)  =0$, we get $A=0$
and $E=-F\tan(\tau)$, so%
\begin{align*}
U_{1}(x)  & =B\sin\left(  \tau x\right)  ,\ \\
U_{3}(x)  & =\frac{F}{\cos\left(  \tau\right)  }\left(  \sin\left(  \tau
x\right)  \cos\left(  \tau\right)  -\cos\left(  \tau x\right)  \sin\left(
\tau\right)  \right)  =\overline{F}\sin\left(  \tau\left(  1-x\right)
\right)  .
\end{align*}
The continuity at $\frac{1}{3}$ and $\frac{2}{3}$ gives%
\[
U_{1}(x)=\frac{C\cos\left(  \frac{\tau}{3}\right)  +D\sin\left(  \frac{\tau
}{3}\right)  }{\sin(\frac{\tau}{3})}\sin\left(  \tau x\right)  ,
\]%
\[
U_{3}(x)=\frac{C\cos\left(  \frac{2\tau}{3}\right)  +D\sin\left(  \frac{2\tau
}{3}\right)  }{\sin(\frac{\tau}{3})}\sin\left(  \tau\left(  1-x\right)
\right)  .
\]
The jump conditions imply that%
\begin{align*}
C  & =\frac{12\sin^{2}\left(  \frac{\tau}{3}\right)  }{\tau-12\sin\left(
\frac{\tau}{3}\right)  \cos\left(  \frac{\tau}{3}\right)  }D,\\
C  & =\frac{12\sin\left(  \frac{\tau}{3}\right)  \sin\left(  \frac{2\tau}%
{3}\right)  -\tau\sin\left(  \tau\right)  }{\tau\cos\left(  \tau\right)
-12sen\left(  \frac{\tau}{3}\right)  \cos\left(  \frac{2\tau}{3}\right)  }D.
\end{align*}
Then there is a countable number of positive values for $\tau$, denoted by
$\{\tau_{k}\}_{k\geq3}$, which are the zeros of the odd function%
\[
g(\tau)=108\sin\frac{1}{3}\tau-36\sin\tau+12\tau\cos\tau-12\tau\cos\frac{1}%
{3}\tau+\allowbreak\tau^{2}\sin\tau=0.
\]
They give rise to the positive eigenvalues $\lambda_{k}=\tau_{k}^{2}$,
$k\geq3$. The fact that $g\left(  \tau\right)  $ has a countable number of
zeros is guaranteed by Lemma \ref{Basis} below.
\end{proof}

\section{Stability of fixed points}

In this section, we will prove the instability of the fixed points $v_{n}%
^{\pm},\ n\geq2,$ by using the linearization around them.

Let us consider the operators $A:H_{0}^{1}\left(  0,1\right)  \rightarrow
H^{-1}\left(  0,1\right)  ,\ B:D(B)\rightarrow H^{-1}\left(  0,1\right)
,\ L:D(L)\rightarrow H^{-1}\left(  0,1\right)  $ given by%
\begin{align*}
Au  & =-u_{xx},\ Bu=-4n\sum_{k=1}^{n-1}\delta_{\frac{k}{n}}u,\ n\in
\mathbb{N},\\
L  & =A+B.
\end{align*}
We observe that $H_{0}^{1}\left(  0,1\right)  \subset D(B)$, so $D(L)=D(A)\cap
D(B)=H_{0}^{1}\left(  0,1\right)  =D(A).$

\begin{lemma}
The operator $A$ is sectorial in the space $H^{-1}\left(  0,1\right)  .$
\end{lemma}

\begin{proof}
We will apply Theorem 36.6 from \cite{SellYou}. Put $V=H_{0}^{1}\left(
0,1\right)  ,\ H=L^{2}(0,1),\ V^{\prime}=H^{-1}(0,1)$. They are Hilbert spaces
and the embeddings $V\subset H\subset V^{\prime}$ are continuous and dense.
The operator $A$ generates the bilinear form $a\left(  u,v\right)  =\int%
_{0}^{1}u_{x}v_{x}dx$. This form is bounded in the $V$-norm and coercive.
Indeed,%
\[
\left\vert a(u,v)\right\vert \leq\left\Vert u\right\Vert _{H_{0}^{1}%
}\left\Vert v\right\Vert _{H_{0}^{1}},
\]%
\[
a(u,u)=\left\Vert u\right\Vert _{H_{0}^{1}}^{2}.
\]
Hence, all conditions in Theorem 36.6 from \cite{SellYou} are satisfied.
\end{proof}

\bigskip

This operator is one to one with bounded inverse $A^{-1}$ \cite[Theorem
36.5]{SellYou}. Also,$\ V^{\prime}=H^{-1}(0,1)$ is a Hilbert space with the
inner product
\[
\left(  u,v\right)  _{H^{-1}}=\left(  A^{-1}u,A^{-1}v\right)  _{H_{0}^{1}},
\]
where $\left(  u,v\right)  _{H_{0}^{1}}=\int_{0}^{1}u_{x}v_{x}dx$ is the inner
product in $H_{0}^{1}(0,1).$ The norm given by this product is equivalent to
the usual norm in $H^{-1}(0,1)$. The operator $A$ is self-adjoint and
positive, which follows from%
\[
\left\langle Au,u\right\rangle _{\left(  H^{-1},H_{0}^{1}\right)  }=\int%
_{0}^{1}\left(  u_{x}\right)  ^{2}dx=\left\Vert u\right\Vert _{H_{0}^{1}}^{2},
\]%
\[
\left\langle Au,v\right\rangle _{\left(  H^{-1},H_{0}^{1}\right)  }=\int%
_{0}^{1}\left(  u_{x}\right)  \left(  v_{x}\right)  dx=\left\langle
u,Av\right\rangle _{\left(  H_{0}^{1},H^{-1}\right)  }.
\]
Let $A^{\alpha}$, $\alpha\in\mathbb{R}$, be the fractional powers of $A$,
where $D(A^{0})=H^{-1}(0,1)$. We take the corresponding spaces $V^{2\alpha
}=D(A^{\alpha})$, denoting their norms by $\left\Vert
\text{\textperiodcentered}\right\Vert _{V^{2\alpha}}$.

\begin{proposition}
\label{EmbeddingCLambda}The embedding $D(A^{\alpha})\subset C^{\lambda
}([0,1])$ is continuous, where $0<\lambda,\alpha\leq1$ are such that
\begin{equation}
0<\frac{\lambda}{2}+\frac{3}{4}<\alpha.\label{AlfaCond}%
\end{equation}

\end{proposition}

\begin{proof}
From Gagliardo-Nirenberg's inequality (see e.g. \cite[Theorem B.3]{SellYou})
we have that%
\[
\left\Vert u\right\Vert _{C^{\lambda}}\leq C_{1}\left\Vert u\right\Vert
_{H_{0}^{1}}^{\theta}\left\Vert u\right\Vert _{L^{2}}^{1-\theta}%
\text{,\ }\forall u\in D(A),
\]
wkere $0<\lambda+\frac{1}{2}<\theta<1$. Since $\left\Vert u\right\Vert
_{H_{0}^{1}}=\left\Vert Au\right\Vert _{H^{-1}}$, we get%
\[
\left\Vert u\right\Vert _{C^{\lambda}}\leq C_{1}\left\Vert Au\right\Vert
_{H^{-1}}^{\theta}\left\Vert u\right\Vert _{L^{2}}^{1-\theta},
\]
and applying Young's inequality it follows that%
\begin{equation}
\left\Vert u\right\Vert _{C^{\lambda}}\leq\varepsilon^{1-\theta}\left\Vert
Au\right\Vert _{H^{-1}}+C_{2}\varepsilon^{-\theta}\left\Vert u\right\Vert
_{L^{2},}\label{IneqTheta}%
\end{equation}
where $C_{4}$ depends on $\theta$. As the norms in $L^{2}(0,1)$ and
$D(A^{\frac{1}{2}})$ are equivalent (see Theorem 37.7 in \cite{SellYou}) we
have that
\[
\left\Vert u\right\Vert _{C^{\lambda}}\leq\varepsilon^{1-\theta}\left\Vert
Au\right\Vert _{H^{-1}}+C_{3}\varepsilon^{-\theta}\left\Vert A^{\frac{1}{2}%
}u\right\Vert _{H^{-1}}.
\]
From the equality%
\[
u=\frac{1}{\Gamma(\alpha)}\int_{0}^{\infty}t^{\alpha-1}e^{-At}A^{\alpha}udt,
\]
applying (\ref{IneqTheta}) we obtain%
\begin{align*}
\left\Vert u\right\Vert _{C^{\lambda}}  & \leq\frac{1}{\Gamma(\alpha)}\int%
_{0}^{\infty}t^{\alpha-1}\left\Vert e^{-At}A^{\alpha}u\right\Vert
_{C^{\lambda}}dt\\
& \leq\frac{1}{\Gamma(\alpha)}\left(  \int_{0}^{\infty}t^{\alpha-1}%
\varepsilon^{1-\theta}\left\Vert Ae^{-At}A^{\alpha}u\right\Vert _{H^{-1}%
}dt+C_{3}\int_{0}^{\infty}t^{\alpha-1}\varepsilon^{-\theta}\left\Vert
A^{\frac{1}{2}}e^{-At}A^{\alpha}u\right\Vert _{H^{-1}}dt\right)  .
\end{align*}
We\ will estimate separately the integrals over the intervals $\left(
0,1\right)  $ and $\left(  1,\infty\right)  $. By using the well-known
inequality%
\[
\left\Vert A^{r}e^{-At}u\right\Vert _{H^{-1}}\leq M_{r}t^{-r}e^{-at}\left\Vert
u\right\Vert _{H^{-1}}%
\]
for some $M_{r},a>0$ \cite[Theorem 37.5]{SellYou} and taking $\varepsilon
=t^{\frac{1}{2}}$ we have first that%
\begin{align*}
& \int_{0}^{1}t^{\alpha-\frac{1}{2}-\frac{\theta}{2}}\left\Vert Ae^{-At}%
A^{\alpha}u\right\Vert _{H^{-1}}dt+C_{3}\int_{0}^{1}t^{\alpha-1-\frac{\theta
}{2}}\left\Vert A^{\frac{1}{2}}e^{-At}A^{\alpha}u\right\Vert _{H^{-1}}dt\\
& \leq C_{4}\int_{0}^{1}t^{\alpha-\frac{3}{2}-\frac{\theta}{2}}\left\Vert
A^{\alpha}u\right\Vert _{H^{-1}}dt=\frac{C_{4}}{\alpha-\frac{1}{2}%
-\frac{\theta}{2}}\left\Vert A^{\alpha}u\right\Vert _{H^{-1}}.
\end{align*}
We oberve that condition (\ref{AlfaCond}) implies the existence of
$\lambda+\frac{1}{2}<\theta<1$ such that $\alpha-\frac{1}{2}-\frac{\theta}%
{2}>0$. Second, setting $\varepsilon=1$ we get%
\begin{align*}
& \int_{1}^{\infty}t^{\alpha-1}\left\Vert Ae^{-At}A^{\alpha}u\right\Vert
_{H^{-1}}dt+C_{3}\int_{1}^{\infty}t^{\alpha-1}\left\Vert A^{\frac{1}{2}%
}e^{-At}A^{\alpha}u\right\Vert _{H^{-1}}dt\\
& \leq C_{5}\left(  \int_{1}^{\infty}t^{\alpha-2}e^{-at}dt+\int_{1}^{\infty
}t^{\alpha-\frac{3}{2}}e^{-at}dt\right)  \left\Vert A^{\alpha}u\right\Vert
_{H^{-1}}\leq C_{6}\left\Vert A^{\alpha}u\right\Vert _{H^{-1}}.
\end{align*}
Hence,%
\[
\left\Vert u\right\Vert _{C^{\lambda}}\leq C_{7}\left\Vert A^{\alpha
}u\right\Vert _{H^{-1}},\ \forall u\in D(A),
\]
whenever (\ref{AlfaCond}) holds. Hence, the inclusion mapping
$J:D(A)\rightarrow C^{\lambda}$ is uniformly bounded in the $A^{\alpha}$-norm
on $D(A)$. As $D(A)$ is dense in \thinspace$D(A^{\alpha})$ for $\alpha\leq1$,
$J$ has a unique continuous extension $J:D(A^{\alpha})\rightarrow C^{\lambda}%
$, which proves the result.
\end{proof}

\bigskip

\begin{lemma}
\label{BoundedOper1}If $3/4<\alpha<1$, then the operator $BA^{-\alpha}%
:H^{-1}(0,1)\rightarrow H^{-1}(0,1)$ is bounded.
\end{lemma}

\begin{proof}
The operator $A^{-\alpha}:H^{-1}(0,1)\rightarrow D(A^{\alpha})$ is bounded for
$\alpha>0$. If $\alpha\in\left(  3/4,1\right)  $, then Proposition
\ref{EmbeddingCLambda} implies that the embedding $D(A^{\alpha})\subset
C([0,1]) $ is continuous. Hence, there exists $C_{1}>0$ such that%
\begin{equation}
\left\Vert A^{-\alpha}u\right\Vert _{C([0,1])}\leq C_{1}\left\Vert
u\right\Vert _{H^{-1}},\ \forall u\in H^{-1}(0,1).\label{IneqC}%
\end{equation}
Then%
\begin{align*}
\left\Vert BA^{-\alpha}u\right\Vert _{H^{-1}}  & =4n\sup_{\phi\in H_{0}%
^{1}(0,1)}\ \frac{\left\vert \left\langle \sum_{k=1}^{n-1}\delta_{\frac{k}{n}%
}A^{-\alpha}u,\phi\right\rangle _{\left(  H^{-1},H_{0}^{1}\right)
}\right\vert }{\left\Vert \phi\right\Vert _{H_{0}^{1}}}\\
& \leq4n\frac{\sum_{k=1}^{n-1}\left\vert A^{-\alpha}u(\frac{k}{n})\right\vert
\left\vert \phi(\frac{k}{n})\right\vert }{\left\Vert \phi\right\Vert
_{H_{0}^{1}}}\\
& \leq4n(n-1)\frac{\left\Vert A^{-\alpha}u\right\Vert _{C([0,1])}\left\Vert
\phi\right\Vert _{C([0,1])}}{\left\Vert \phi\right\Vert _{H_{0}^{1}}}\leq
C_{2}\left\Vert u\right\Vert _{H^{-1}}.
\end{align*}

\end{proof}

\begin{lemma}
The operator $L$ is sectorial in the space $H^{-1}\left(  0,1\right)  .$
\end{lemma}

\begin{proof}
Since $BA^{-\alpha}:H^{-1}(0,1)\rightarrow H^{-1}(0,1)$ is bounded for
$\alpha\in\left(  3/4,1\right)  $ and the eigenvalues of $A$ are positive, the
result follows from Corollary 1.4.5 in \cite{Henry}.
\end{proof}

\begin{lemma}
\label{AcotB}For any $\varepsilon>0$ there is $K\left(  \varepsilon\right)
>0$ such that%
\[
\left\Vert Bu\right\Vert _{H^{-1}}\leq\varepsilon\left\Vert u_{x}\right\Vert
_{L^{2}}+K(\varepsilon)\left\Vert u\right\Vert _{H^{-1}}\text{ }\forall u\in
H_{0}^{1}(0,1).
\]

\end{lemma}

\begin{proof}
The element $u\in H^{-1}(0,1)$ is represented by $u=df/dx$ with $f\in
L^{2}(0,1)$. Then%
\begin{align*}
\left\Vert Bu\right\Vert _{H^{-1}}  & =\left\Vert 4n\sum_{k=1}^{n-1}%
\delta_{\frac{k}{n}}u\right\Vert _{H^{-1}}=4n\sup_{\phi\in H_{0}^{1}%
(0,1)}\ \frac{\left\vert \sum_{k=1}^{n-1}\left\langle \delta_{\frac{k}{n}%
}u,\phi\right\rangle _{\left(  H^{-1},H_{0}^{1}\right)  }\right\vert
}{\left\Vert \phi\right\Vert _{H_{0}^{1}}}\\
& \leq4n\sup_{\phi\in H_{0}^{1}(0,1)}\frac{\sum_{k=1}^{n-1}\left\vert
u(\frac{k}{n})\right\vert \left\vert \phi(\frac{k}{n})\right\vert }{\left\Vert
\phi\right\Vert _{H_{0}^{1}}}\leq4n\sum_{k=1}^{n-1}\left\vert u(\frac{k}%
{n})\right\vert ,
\end{align*}
where the last inequality follows from $\left\vert \phi(\frac{k}%
{n})\right\vert =\left\vert \int_{0}^{\frac{k}{n}}\phi_{x}dx\right\vert
\leq\left\Vert \phi\right\Vert _{H_{0}^{1}}.$ We recall the
Gagliardo-Nirenberg inequality in dimension $1$:%
\[
\left\Vert \frac{d^{j}f}{dx^{j}}\right\Vert _{L^{p}}\leq C_{1}\left\Vert
\frac{d^{m}f}{dx^{m}}\right\Vert _{L^{r}}^{a}\left\Vert f\right\Vert _{L^{q}%
}^{1-a},
\]
where $0\leq j<m$, $1\leq q,r,p\leq\infty$, $1/p=j+a(1/r-m)+(1-a)/q$. Choosing
$m=r=q=2,\ p=\infty,\ a=3/4,\ j=1$ we have%
\[
\left\vert u(\frac{k}{n})\right\vert =\left\vert \frac{df}{dx}(\frac{k}%
{n})\right\vert \leq\left\Vert \frac{df}{dx}\right\Vert _{L^{\infty}}\leq
C_{1}\left\Vert f\right\Vert _{H^{2}}^{\frac{3}{4}}\left\Vert f\right\Vert
_{L^{2}}^{\frac{1}{4}}.
\]
Further, we need to estimate the norm $\left\Vert f\right\Vert _{H^{2}}$. We
have:%
\begin{align*}
\left\Vert f\right\Vert _{H^{2}}  & =\left\Vert f\right\Vert _{L^{2}%
}+\left\Vert f_{x}\right\Vert _{L^{2}}+\left\Vert f_{xx}\right\Vert _{L^{2}}\\
& =\left\Vert u\right\Vert _{H^{-1}}+\left\Vert u\right\Vert _{L^{2}%
}+\left\Vert u_{x}\right\Vert _{L^{2}}\\
& \leq\left\Vert u\right\Vert _{H^{-1}}+\left(  1+\frac{1}{\sqrt{\lambda_{1}}%
}\right)  \left\Vert u_{x}\right\Vert _{L^{2}},
\end{align*}
and then $\left\Vert u\right\Vert _{H^{-1}}=\left\Vert f\right\Vert _{L^{2}}$
gives for an arbitrarily small $\varepsilon>0$ that%
\begin{align*}
\left\Vert Bu\right\Vert _{H^{-1}}  & \leq C_{2}\left(  \left\Vert
u\right\Vert _{H^{-1}}+\left\Vert u_{x}\right\Vert _{L^{2}}^{\frac{3}{4}%
}\left\Vert u\right\Vert _{H^{-1}}^{\frac{1}{4}}\right) \\
& \leq\varepsilon\left\Vert u_{x}\right\Vert _{L^{2}}+K(\varepsilon)\left\Vert
u\right\Vert _{H^{-1}}.
\end{align*}

\end{proof}

\begin{lemma}
\label{Basis}The operator $L$ has a discrete spectrum $\lambda_{1}\leq
\lambda_{2}\leq...\leq\lambda_{k}\leq...$ with $\lambda_{k}\rightarrow+\infty
$. The eigenfunctions $\{e_{1},e_{2},...\}$ of the operator $L$ form an
orthonormal basis of $H^{-1}(0,1).$
\end{lemma}

\begin{proof}
For $\alpha>0$ we define the operator $\widetilde{L}u=Lu+\alpha u$. We
consider now the problem%
\begin{equation}
\widetilde{L}u=h\in H^{-1}(0,1).\label{EqL}%
\end{equation}
Let $\widetilde{a}:H_{0}^{1}(0,1)\times H_{0}^{1}(0,1)\rightarrow\mathbb{R}$
be the bilinear form given by $\widetilde{a}(u,v)=\int_{0}^{1}u_{x}%
v_{x}dx-\left\langle Bu,v\right\rangle _{\left(  H^{-1},H_{0}^{1}\right)
}+\alpha\int_{0}^{1}uvdx$. This form is continuous and coercive for $\alpha$
great enough. Indeed, using Lemma \ref{AcotB} we get%
\[
\left\vert \widetilde{a}(u,v)\right\vert \leq\left\Vert u\right\Vert
_{H_{0}^{1}}\left\Vert v\right\Vert _{H_{0}^{1}}+4n\sum_{k=1}^{n-1}\left\vert
u(\frac{k}{n})\right\vert \left\vert v(\frac{k}{n})\right\vert +\alpha
\left\Vert u\right\Vert _{L^{2}}\left\Vert v\right\Vert _{L^{2}}\leq
C_{1}\left\Vert u\right\Vert _{H_{0}^{1}}\left\Vert v\right\Vert _{H_{0}^{1}},
\]%
\begin{align*}
\widetilde{a}(u,u)  & \geq\left\Vert u\right\Vert _{H_{0}^{1}}^{2}-\left\Vert
Bu\right\Vert _{H^{-1}}\left\Vert u\right\Vert _{H_{0}^{1}}+\alpha\left\Vert
u\right\Vert _{L^{2}}^{2}\\
& =\left\Vert u\right\Vert _{H_{0}^{1}}^{2}-\varepsilon\left\Vert u\right\Vert
_{H_{0}^{1}}^{2}-K(\varepsilon)\left\Vert u\right\Vert _{H^{-1}}\left\Vert
u\right\Vert _{H_{0}^{1}}+\alpha\left\Vert u\right\Vert _{L^{2}}^{2}\\
& \geq\left\Vert u\right\Vert _{H_{0}^{1}}^{2}-2\varepsilon\left\Vert
u\right\Vert _{H_{0}^{1}}^{2}+\left(  \alpha-K_{1}(\varepsilon)\right)
\left\Vert u\right\Vert _{L^{2}}^{2}.
\end{align*}
Taking $\varepsilon=1/4$ and $\alpha>K_{1}(\varepsilon)$ we obtain the result.
Hence, Lax-Milgram's theorem implies the existence of a unique solution to
problem (\ref{EqL}).

Moreover,%
\[
(1-2\varepsilon)\left\Vert u\right\Vert _{H_{0}^{1}}^{2}\leq\widetilde{a}%
(u,u)=\left\langle \widetilde{L}u,u\right\rangle _{\left(  H^{-1},H_{0}%
^{1}\right)  }\leq\left\Vert h\right\Vert _{H^{-1}}\left\Vert u\right\Vert
_{H_{0}^{1}}\leq\frac{1}{2(1-2\varepsilon)}\left\Vert h\right\Vert _{H^{-1}%
}+\frac{1-2\varepsilon}{2}\left\Vert u\right\Vert _{H_{0}^{1}}^{2},
\]
so it follows that $L^{-1}:H^{-1}(0,1)\rightarrow H_{0}^{1}(0,1)$ is bounded.
The compact embedding $H_{0}^{1}(0,1)\subset H^{-1}(0,1)$ gives that
$\widetilde{L}^{-1}:H^{-1}(0,1)\rightarrow H^{-1}(0,1)$ is compact. Thus,
$\widetilde{L}$ possesses a discrete spectrum $\widetilde{\lambda}_{1}%
\leq\widetilde{\lambda}_{2}\leq...$ such that $\widetilde{\lambda}%
_{n}\rightarrow\infty$ and the eigenfunctions $\{e_{1},e_{2},...\}$ form and
orthonormal basis of $H^{-1}(0,1)$ \cite[Corollary 3.26]{Robinson}.

The result follows from the fact that the eigenvalues of the operator $L$ are
$\lambda_{k}=\widetilde{\lambda}_{k}-\alpha$ and the eigenfunctions coincide.
\end{proof}

\bigskip

\begin{lemma}
\label{EigenvaluesL}$L$ has at least one negative eigenvalue. The number of
negative eigenvalues is finite and $0$ is not an eigenvalue.
\end{lemma}

\begin{proof}
In Lemmas \ref{EigenvaluesU2}, \ref{EigenvaluesU3} these results were
established for $n=2,3$. Let us prove them for any $n\geq2.$

From Lemma \ref{Basis} we know that $L$ has a discrete spectrum $\lambda
_{1}\leq\lambda_{2}\leq...\leq\lambda_{k}\leq...$ with $\lambda_{k}%
\rightarrow+\infty$, from where we deduce immediately that the number of
negative eigenvalues is at most finite.

By contradiction, assume that $\lambda_{1}\geq0$. Then, since the
eigenfunctions form a basis in $H^{-1}(0,1)$, it follows easily that
$\left\langle Lu,u\right\rangle _{\left(  H^{-1},H_{0}^{1}\right)  }\geq0$ for
any $u\in H_{0}^{1}\left(  0,1\right)  .$ Let%
\[
u\left(  x\right)  =\left\{
\begin{array}
[c]{c}%
nx\text{ if }0\leq x\leq\frac{1}{n},\\
2-nx\text{ if }\frac{1}{n}\leq x\leq\frac{2}{n},\\
0\text{ if }\frac{2}{n}\leq x\leq1.
\end{array}
\right.
\]
Then $u\in H_{0}^{1}\left(  0,1\right)  $ and $u_{xx}=-2n\delta_{\frac{1}{n}%
}+n\delta_{\frac{2}{n}}$. Thus,%
\begin{align*}
\left\langle Lu,u\right\rangle _{\left(  H^{-1},H_{0}^{1}\right)  }  &
=\left\langle 2n\delta_{\frac{1}{n}}-n\delta_{\frac{2}{n}},u\right\rangle
_{\left(  H^{-1},H_{0}^{1}\right)  }-4n\left\langle \sum_{k=1}^{n-1}%
\delta_{\frac{k}{n}}u,u\right\rangle _{\left(  H^{-1},H_{0}^{1}\right)  }\\
& =2nu\left(  \frac{1}{n}\right)  -4nu^{2}\left(  \frac{1}{n}\right)  =-2n<0,
\end{align*}
which is a contradiction, showing that $\lambda_{1}<0$.

Finally, we will prove that $0$ is not an eigenvalue. From (\ref{taon}) an
eigenfunction $U\left(  x\right)  $ correspoding to $\lambda=0$ has to satisfy%
\[
U\left(  x\right)  =U_{k}\left(  x\right)  =A_{k}+B_{k}x\text{ if }x\in\left[
\frac{k-1}{n},\frac{k}{n}\right]  \text{, }k=1,...,n.
\]
From $U(0)=0$ we have $A_{1}=0$ and by $U(1)=0$ and the continuity and jump
conditions we obtain the following equations:%
\begin{align*}
A_{k}+\frac{k}{n}B_{k}  & =A_{k+1}+\frac{k}{n}B_{k+1},\\
B_{k}-B_{k+1}  & =4n\left(  A_{k+1}+\frac{k}{n}B_{k+1}\right)
,\ k=1,...,n-1,\\
A_{n}+B_{n}  & =0.
\end{align*}
Then we get the linear system $Mz=0$, where $z=(B_{1},A_{2},B_{2}%
,...,A_{n},B_{n})^{t}$ and
\[
M=\left(
\begin{array}
[c]{ccccccccccc}%
\frac{1}{n} & -1 & -\frac{1}{n} & 0 & \cdots & \cdots & \cdots & \cdots &
\cdots & \cdots & 0\\
1 & -4n & -5 & 0 & 0 & \cdots & \cdots & \cdots & \cdots & \cdots & 0\\
0 & 1 & \frac{2}{n} & -1 & -\frac{2}{n} & 0 & \cdots & \cdots & \cdots &
\cdots & 0\\
\vdots & 0 & 1 & -4n & -9 & 0 & 0 & \cdots & \cdots & \cdots & 0\\
\vdots & \vdots & 0 & \ddots & \ddots & \ddots & \ddots & \ddots & \cdots &
\cdots & \vdots\\
\vdots & \vdots & \vdots & 0 & 1 & \frac{k}{n} & -1 & -\frac{k}{n} & 0 &
\cdots & \vdots\\
\vdots & \vdots & \vdots & \vdots & 0 & 1 & -4n & -1-4k & 0 & \ddots &
\vdots\\
\vdots & \vdots & \vdots & \vdots & \vdots & \ddots & \ddots & \ddots & \ddots
& \ddots & 0\\
\vdots & \vdots & \vdots & \vdots & \vdots & \vdots & \ddots & 1 & \frac
{n-1}{n} & -1 & -\frac{n-1}{n}\\
\vdots & \vdots & \vdots & \vdots & \vdots & \vdots & \vdots & 0 & 1 & -4n &
-1-4(n-1)\\
0 & 0 & 0 & 0 & \cdots & \cdots & \cdots & \cdots & 0 & 1 & 1
\end{array}
\right)  .
\]
Applying Gauss' method we transform $M$ into\ the following superior
triangular matrix:%
\[
\widetilde{M}=\left(
\begin{array}
[c]{ccccccccccc}%
\frac{1}{n} & -1 & -\frac{1}{n} & 0 & \cdots & \cdots & \cdots & \cdots &
\cdots & \cdots & 0\\
0 & -3n & -4 & 0 & 0 & \cdots & \cdots & \cdots & \cdots & \cdots & 0\\
0 & 0 & \frac{2}{3n} & -1 & -\frac{2}{n} & 0 & \cdots & \cdots & \cdots &
\cdots & 0\\
\vdots & 0 & 0 & -\frac{5n}{2} & -6 & 0 & 0 & \cdots & \cdots & \cdots & 0\\
\vdots & \vdots & 0 & \ddots & \ddots & \ddots & \ddots & \ddots & \cdots &
\cdots & \vdots\\
\vdots & \vdots & \vdots & 0 & 0 & \frac{k}{(2k-1)n} & -1 & -\frac{k}{n} & 0 &
\cdots & \vdots\\
\vdots & \vdots & \vdots & \vdots & 0 & 0 & -\frac{1+2k}{k}n & -2-2k & 0 &
\ddots & \vdots\\
\vdots & \vdots & \vdots & \vdots & \vdots & \ddots & \ddots & \ddots & \ddots
& \ddots & 0\\
\vdots & \vdots & \vdots & \vdots & \vdots & \vdots & \ddots & 0 & \frac
{n-1}{(2n-3)n} & -1 & -\frac{n-1}{n}\\
\vdots & \vdots & \vdots & \vdots & \vdots & \vdots & \vdots & 0 & 0 &
-\frac{2n-1}{n-1}n & -2n\\
0 & 0 & 0 & 0 & \cdots & \cdots & \cdots & \cdots & 0 & 0 & \frac{1}{2n-1}%
\end{array}
\right)  .
\]
Let us prove this fact by induction. Assume that the result is true for the
row $2k$. Then, applying Gauss' method we obtain the following non-zero
elements in the rows $2k+1$ and $2k+2:$%
\begin{align*}
\widetilde{M}_{2k+1,2k+1}  & =\frac{k+1}{n}-\frac{k\left(  2+2k\right)
}{(1+2k)n}=\frac{k+1}{(2(k+1)-1)n},\\
\widetilde{M}_{2k+1,2k+2}  & =-1,\ \widetilde{M}_{2k+1,2k+3}=-\frac{k+1}{n},\\
\widetilde{M}_{2k+2,2k+2}  & =-4n+\frac{\left(  2(k+1)-1\right)  n}%
{k+1}=-\frac{1+2(k+1)}{k+1}n,\\
\widetilde{M}_{2k+2,2k+3}  & =-1-4(k+1)+\frac{(2(k+1)-1)n}{k+1}\frac{k+1}%
{n}=-2-2(k+1).
\end{align*}
Since $\det(\widetilde{M})\not =0$, the only solution is $z=0$, so $\lambda=0$
is not an eigenvalue.
\end{proof}

\bigskip

We know from Proposition \ref{EigenvaluesL} that $L$ has a discrete spectrum
$\delta(L)$ given by the values%
\[
\lambda_{1}\leq\lambda_{2}\leq...\leq\lambda_{m}<0<\lambda_{m+1}\leq
\lambda_{m+2}<...
\]
We denote $\sigma_{1}=\{\lambda_{1},...,\lambda_{m}\}$, $\sigma_{2}%
=\sigma(L)\backslash\sigma_{1}$. Let $E_{1},E_{2}$ be the projections
associated to these spectral sets and $V_{j}^{\prime}=E_{j}(H^{-1}(0,1))$.
Then \cite[Theorems 1.5.2, 1.5.3]{Henry} $H^{-1}(0,1)=V_{1}^{\prime}\oplus
V_{2}^{\prime}$, the spaces $V_{j}^{\prime}$ are invariant under the operator
$L$ and if $L_{j}$ are the restrictions of $L$ to $V_{j}^{\prime}$, then
\[
L_{1}:V_{1}^{\prime}\rightarrow V_{1}^{\prime}\text{ is bounded, }\sigma
(L_{1})=\sigma_{1},
\]%
\[
D(L_{2})=D(L)\cap V_{2}^{\prime},\ \sigma(L_{2})=\sigma_{2},
\]%
\[
L_{2}\text{ is sectorial,}%
\]%
\begin{equation}
\left\Vert e^{-L_{1}t}\right\Vert _{\mathcal{L}(H^{-1})}\leq Ce^{-\lambda
_{m}t}\ \forall t\leq0,\label{IneqExpL1}%
\end{equation}%
\begin{equation}
\left\Vert e^{-L_{2}t}\right\Vert _{\mathcal{L}(H^{-1})}\leq Ce^{-\lambda
_{m+1}t}\ \forall t>0,\label{IneqExpL2A}%
\end{equation}%
\begin{equation}
\left\Vert L_{2}e^{-L_{2}t}\right\Vert _{\mathcal{L}(H^{-1})}\leq
Ct^{-1}e^{-\lambda_{m+1}t}\ \forall t>0,\label{IneqExpL2B}%
\end{equation}
for some constant $C>0.$ As $L_{2}$ is sectorial, it is the infinitesimal
generator of the analitic strongly continuous semigroup of bounded linear
operators $e^{-L_{2}t}$ \cite[Theorem 1.3.4]{Henry}. The same is valid for the
operator $L_{1}$, because being $L_{1}$ a linear bounded operator, it is
sectorial as well \cite[p. 19]{Henry}.

\begin{lemma}
\label{IneqL2}For any $\alpha\in\lbrack0,1)$ the following inequalities are
satisfied:%
\begin{equation}
\left\Vert e^{-L_{2}t}u\right\Vert _{V^{2\alpha}}\leq Ct^{-\alpha}%
e^{-\lambda_{m+1}t}\left\Vert u\right\Vert _{H^{-1}},\ \forall u\in
V_{2}^{\prime},\ t>0,\label{IneqL2A}%
\end{equation}%
\begin{equation}
\left\Vert e^{-L_{2}t}u\right\Vert _{V^{2\alpha}}\leq Ce^{-\lambda_{m+1}%
t}\left\Vert u\right\Vert _{V^{2\alpha}},\ \forall u\in V_{2}^{\prime}\cap
D(A^{\alpha}),\ t\geq0,\label{IneqL2B}%
\end{equation}%
\begin{equation}
\left\Vert e^{-L_{1}t}u\right\Vert _{V^{2\alpha}}\leq Ce^{-\lambda_{m}%
t}\left\Vert u\right\Vert _{H^{-1}},\ \forall u\in V_{1}^{\prime}%
,\ t\leq0,\label{IneqL2C}%
\end{equation}%
\begin{equation}
\left\Vert e^{-L_{1}t}u\right\Vert _{V^{2\alpha}}\leq Ce^{-\lambda_{m}%
t}\left\Vert u\right\Vert _{V^{2\alpha}},\ \forall u\in V_{1}^{\prime}%
,\ t\leq0,\label{IneqL2D}%
\end{equation}
where $C>0$ is some constant.
\end{lemma}

\begin{proof}
As $D(A)=H_{0}^{1}(0,1)=D(L)$, the eigenvalues of $A$ are positive and the
operator $\left(  A-L\right)  A^{-\alpha}=BA^{-\alpha}$ is bounded by Lemma
\ref{BoundedOper1}, inequalities (\ref{IneqL2A})-(\ref{IneqL2B}) follow from
Theorem 1.5.4 in \cite{Henry}. Inequalities (\ref{IneqL2C}), (\ref{IneqL2D})
are a consequence of (\ref{IneqExpL1}) and the fact that the norms $\left\Vert
\text{\textperiodcentered}\right\Vert _{H^{-1}},\ \left\Vert
\text{\textperiodcentered}\right\Vert _{V^{2\alpha}}$, $0\leq\alpha\leq1 $,
are equivalent in $V_{1}^{\prime}$.
\end{proof}

\bigskip

\begin{lemma}
\label{DisL1L2}The operators $L_{1},L_{2}$ satisfy%
\begin{equation}
\left(  L_{1}v,v\right)  _{H^{-1}}\geq\lambda_{1}\left\Vert v\right\Vert
_{H^{-1}}\ \forall v\in V_{1}^{\prime},\label{L1Ineq0}%
\end{equation}%
\begin{equation}
\left(  L_{1}v,v\right)  _{H^{-1}}\leq\lambda_{m}\left\Vert v\right\Vert
_{H^{-1}}\ \forall v\in V_{1}^{\prime},\label{L1Eq}%
\end{equation}%
\begin{equation}
\left(  L_{2}v,v\right)  _{H^{-1}}\geq\lambda_{m+1}\left\Vert v\right\Vert
_{H^{-1}}\ \forall v\in D(L_{2}).\label{L2Ineq}%
\end{equation}

\end{lemma}

\begin{proof}
(\ref{L1Eq}) is a consequence of $L_{1}v=\lambda_{1}v$ for any $v\in
V_{1}^{\prime}$.

In view of Lemma \ref{Basis}, the subsets of eigenfunctions $\{e_{1}%
,...,e_{m}\}$ and$\{e_{m+1},e_{m+2},...\}$ form an orthonormal basis in
$V_{1}^{\prime}$ and $V_{2}^{\prime}$, respectively. Then
\[
\left(  L_{1}v,v\right)  _{H^{-1}}=\sum_{i=1}^{m}\lambda_{i}\left(
v,e_{i}\right)  _{H^{-1}}^{2}\leq\lambda_{m}\sum_{i=m+1}^{\infty}\left(
v,e_{i}\right)  _{H^{-1}}^{2}=\lambda_{m}\left\Vert v\right\Vert _{H^{-1}%
},\ \forall v\in V_{1}^{\prime},
\]%
\[
\left(  L_{1}v,v\right)  _{H^{-1}}=\sum_{i=1}^{m}\lambda_{i}\left(
v,e_{i}\right)  _{H^{-1}}^{2}\geq\lambda_{1}\sum_{i=m+1}^{\infty}\left(
v,e_{i}\right)  _{H^{-1}}^{2}=\lambda_{1}\left\Vert v\right\Vert _{H^{-1}%
},\ \forall v\in V_{1}^{\prime},
\]%
\[
\left(  L_{2}v,v\right)  _{H^{-1}}=\sum_{i=m+1}^{\infty}\lambda_{i}\left(
v,e_{i}\right)  _{H^{-1}}^{2}\geq\lambda_{m+1}\sum_{i=m+1}^{\infty}\left(
v,e_{i}\right)  _{H^{-1}}^{2}=\lambda_{m+1}\left\Vert v\right\Vert _{H^{-1}%
}\text{, }\forall v\in D(L_{2}).
\]

\end{proof}

\bigskip

We define $A_{H}:D(A_{H})\rightarrow L^{2}\left(  0,1\right)  $ as the
operator $A_{H}u=-u_{xx}$ with $D(A)=H^{2}\left(  0,1\right)  \cap H_{0}%
^{1}\left(  0,1\right)  $. Denote $V_{H}^{2r}=D\left(  A_{H}^{r}\right)  $,
where $A_{H}^{r}$,. $r\in\mathbb{R}$, are the fractional powers of the
operator $A_{H}$. We take $3/4<r<1$. Since for such $r$ we have that
$V_{H}^{2r}\subset C^{1}([0,1])$, one can choose a neighborhood $O_{\delta
}\left(  v_{n}^{\pm}\right)  $ in $V_{H}^{2r}$ such that any $v\in O_{\delta
}\left(  v_{n}^{+}\right)  $ ($O_{\delta}\left(  v_{n}^{-}\right)  $) satisfies:

\begin{itemize}
\item $v$ has exactly $n+1$ zeros $x_{v}^{j},$ $j=0,...,n,$ in $[0,1]$, where
$x_{v}^{0}=0,\ x_{v}^{n}=1.$

\item $v\left(  x\right)  >0$ ($<0$), $\forall x\in\left(  x_{v}^{j}%
,x_{v}^{j+1}\right)  $, if $j$ is even, and $v\left(  x\right)  $ $<0$ ($>0$),
$\forall x\in\left(  x_{v}^{j},x_{v}^{j+1}\right)  $, if $j$ is odd, for
$j=0,...,n-1$.
\end{itemize}

Also, $x_{v_{k}}^{j}\rightarrow\frac{j}{n}$ if $v_{k}\rightarrow v_{n}^{+}$ in
$V^{2r}.$ Moreover, this neighborhood can be chosen such that the solutions to
problem (\ref{Eq}) starting at $O_{\delta}\left(  v_{n}^{\pm}\right)  $ are
unique while they remain inside $O_{\delta}\left(  v_{n}^{\pm}\right)  $
\cite[p.35]{MorVal}. For $v\in O_{\delta}\left(  v_{n}^{\pm}\right)  $ let
$z=v-v_{n}^{\pm}$.

\begin{lemma}
\label{AcotzLemma2}Let $n\geq2$. There is $\rho\leq\delta$ such that if $v\in
O_{\rho}(v_{n}^{\pm})$, then%
\begin{equation}
\left\Vert H_{0}\left(  v\left(  \text{\textperiodcentered}\right)  \right)
-H_{0}\left(  v_{n}^{\pm}\left(  \text{\textperiodcentered}\right)  \right)
-4n\sum_{j=1}^{n-1}\delta_{\frac{j}{n}}z\left(  \text{\textperiodcentered
}\right)  \right\Vert _{H^{-1}}\leq C\left\Vert z\right\Vert _{C^{1}%
([0,1])}^{\frac{1}{2}}\left\Vert z\right\Vert _{C([0,1])},\label{AcotzGen}%
\end{equation}
for some positive contant $C.$
\end{lemma}

\begin{proof}
Let us consider $v_{n}^{+}$. It is easy to see that%
\[
H_{0}\left(  v_{n}^{+}\left(  x\right)  \right)  =\left\{
\begin{array}
[c]{c}%
1\text{ for }x\in\left(  \frac{j}{n},\frac{j+1}{n}\right)  ,\text{ if }j\text{
is even,}\\
-1\text{ for }x\in\left(  \frac{j}{n},\frac{j+1}{n}\right)  ,\text{ if
}j\text{ is odd,}%
\end{array}
\right.
\]%
\[
H_{0}\left(  v\left(  x\right)  \right)  =\left\{
\begin{array}
[c]{c}%
1\text{ for }x\in\left(  x_{v}^{j},x_{v}^{j+1}\right)  ,\ \text{if }j\text{ is
even,}\\
-1\text{ for }x\in\left(  x_{v}^{j},x_{v}^{j+1}\right)  ,\ \text{if }j\text{
is odd,}%
\end{array}
\right.
\]
where $j=0,...,n-1$. Then%
\[
H_{0}\left(  v\left(  x\right)  \right)  -H_{0}\left(  v_{n}^{+}\left(
x\right)  \right)  =\left\{
\begin{array}
[c]{c}%
2(-1)^{j+1}\text{ for }x\in\left(  \frac{j}{n},x_{v}^{j}\right)  \text{, if
}x_{v}^{j}>\frac{j}{n}\text{,}\\
-2(-1)^{j+1}\text{ for }x\in\left(  x_{v}^{j},\frac{j}{n}\right)  \text{, if
}x_{v}^{j}<\frac{j}{n}\text{,}\\
0\text{, otherwise,}%
\end{array}
\right.
\]
and for $\varphi\in H_{0}^{1}\left(  0,1\right)  ,$%
\[
\left\langle H_{0}\left(  v\left(  \text{\textperiodcentered}\right)  \right)
-H_{0}\left(  v_{n}^{+}\left(  \text{\textperiodcentered}\right)  \right)
,\varphi\right\rangle _{\left(  H^{-1},H_{0}^{1}\right)  }=\sum_{j=1}%
^{n-1}\int_{\frac{j}{n}}^{x_{v}^{j}}2(-1)^{j+1}\varphi\left(  x\right)
dx=2\sum_{j=1}^{n-1}\left(  x_{v}^{j}-\frac{j}{n}\right)  (-1)^{j+1}%
\varphi\left(  x_{v,j}^{\ast}\right)  ,
\]
for some $x_{v,j}^{\ast}\in\lbrack\min\{x_{v}^{j},\frac{j}{n}\},\max
\{x_{v}^{j},\frac{j}{n}\}]$. Hence,%
\begin{align*}
& \left\langle H_{0}\left(  v\left(  \text{\textperiodcentered}\right)
\right)  -H_{0}\left(  v_{n}^{+}\left(  \text{\textperiodcentered}\right)
\right)  -4n\sum_{j=1}^{n-1}\delta_{\frac{j}{n}}\left(  v\left(
\text{\textperiodcentered}\right)  -v_{n}^{+}\left(  \text{\textperiodcentered
}\right)  \right)  ,\varphi\right\rangle _{\left(  H^{-1},H_{0}^{1}\right)
}\\
& =\sum_{j=1}^{n-1}\left(  2\left(  x_{v}^{j}-\frac{j}{n}\right)
(-1)^{j+1}\varphi\left(  x_{v,j}^{\ast}\right)  -4nv\left(  \frac{j}%
{n}\right)  \varphi\left(  \frac{j}{n}\right)  \right)  .
\end{align*}
On the other hand,%
\begin{equation}
-v\left(  \frac{j}{n}\right)  =v\left(  x_{v}^{j}\right)  -v\left(  \frac
{j}{n}\right)  =v^{\prime}\left(  \widetilde{x}_{v,j}\right)  \left(
x_{v}^{j}-\frac{j}{n}\right)  ,\label{Eqvj}%
\end{equation}
for some $\widetilde{x}_{v,j}\in\lbrack\min\{x_{v}^{j},\frac{j}{n}%
\},\max\{x_{v}^{j},\frac{j}{n}\}]$, so%
\begin{align}
& \left\vert \left\langle H_{0}\left(  v\left(  \text{\textperiodcentered
}\right)  \right)  -H_{0}\left(  v_{n}^{+}\left(  \text{\textperiodcentered
}\right)  \right)  -4n\sum_{j=1}^{n-1}\delta_{\frac{j}{n}}\left(  v\left(
\text{\textperiodcentered}\right)  -v_{n}^{+}\left(  \text{\textperiodcentered
}\right)  \right)  ,\varphi\right\rangle _{\left(  H^{-1},H_{0}^{1}\right)
}\right\vert \nonumber\\
& \leq2\sum_{j=1}^{n-1}\left\vert (-1)^{j+1}\frac{\varphi\left(  x_{v,j}%
^{\ast}\right)  }{v^{\prime}\left(  \widetilde{x}_{v,j}\right)  }%
+2n\varphi\left(  \frac{j}{n}\right)  \right\vert \left\vert v\left(  \frac
{j}{n}\right)  -v_{n}^{+}\left(  \frac{j}{n}\right)  \right\vert \nonumber\\
& =2\sum_{j=1}^{n-1}\left\vert (-1)^{j+1}\frac{\varphi\left(  \frac{j}%
{n}\right)  +\int_{\frac{j}{n}}^{x_{v,j}^{\ast}}\varphi^{\prime}\left(
x\right)  dx}{v^{\prime}\left(  \widetilde{x}_{v,j}\right)  }+2n\varphi\left(
\frac{j}{n}\right)  \right\vert \left\vert v\left(  \frac{j}{n}\right)
-v_{n}^{+}\left(  \frac{j}{n}\right)  \right\vert \nonumber\\
& \leq C_{1}\sum_{j=1}^{n-1}\left(  \left\vert \frac{(-1)^{j+1}}{v^{\prime
}\left(  \widetilde{x}_{v,j}\right)  }+2n\right\vert +\frac{\left\vert
x_{v,j}-\frac{j}{n}\right\vert ^{\frac{1}{2}}}{\left\vert v^{\prime}\left(
\widetilde{x}_{v,j}\right)  \right\vert }\right)  \left\Vert v-v_{n}%
^{+}\right\Vert _{C([0,1])}\left\Vert \varphi\right\Vert _{H_{0}^{1}%
}\nonumber\\
& \leq C_{1}\sum_{j=1}^{n-1}\left(  \frac{2n}{\left\vert v^{\prime}\left(
\widetilde{x}_{v,j}\right)  \right\vert }\left(  \left\vert \left(  v_{n}%
^{+}\right)  ^{\prime}\left(  \widetilde{x}_{v,j}\right)  -\left(  v_{n}%
^{+}\right)  ^{\prime}\left(  \frac{j}{n}\right)  \right\vert +\left\vert
v^{\prime}\left(  \widetilde{x}_{v,j}\right)  -\left(  v_{n}^{+}\right)
^{\prime}\left(  \widetilde{x}_{v,j}\right)  \right\vert \right)
+\frac{\left\Vert z\right\Vert _{C([0,1])}^{\frac{1}{2}}}{\left\vert
v^{\prime}\left(  \widetilde{x}_{v,j}\right)  \right\vert ^{\frac{3}{2}}%
}\right)  \left\Vert z\right\Vert _{C([0,1])}\left\Vert \varphi\right\Vert
_{H_{0}^{1}}\nonumber\\
& \leq C_{2}\sum_{j=1}^{n-1}\left(  \frac{2n\left\Vert v_{n}^{+}\right\Vert
_{H^{2}}\left\Vert z\right\Vert _{C([0,1])}^{\frac{1}{2}}}{\left\vert
v^{\prime}\left(  \widetilde{x}_{v,j}\right)  \right\vert ^{\frac{3}{2}}%
}+\frac{2n\left\Vert z\right\Vert _{C^{1}([0,1])}}{\left\vert v^{\prime
}\left(  \widetilde{x}_{v,j}\right)  \right\vert }+\frac{\left\Vert
z\right\Vert _{C([0,1])}^{\frac{1}{2}}}{\left\vert v^{\prime}\left(
\widetilde{x}_{v,j}\right)  \right\vert ^{\frac{3}{2}}}\right)  \left\Vert
z\right\Vert _{C([0,1])}\left\Vert \varphi\right\Vert _{H_{0}^{1}%
},\label{IneqGen}%
\end{align}
where we have used that $\left(  v_{n}^{+}\right)  ^{\prime}\left(  \frac
{j}{n}\right)  =\frac{\left(  -1\right)  ^{j}}{2n}$. The embedding $V_{H}%
^{2r}\subset C^{1}([0,1])$ implies the existence of $\gamma_{0}>0$ and
$\rho\leq\delta$ such that if $\left\Vert v-u_{n}^{+}\right\Vert _{V^{2r}}%
\leq\rho$, then $\left\vert v^{\prime}\left(  \widetilde{x}_{v,j}\right)
\right\vert \geq\gamma_{0}$ for all $j$. Hence,%
\begin{align*}
& \left\vert \left\langle H_{0}\left(  v\left(  \text{\textperiodcentered
}\right)  \right)  -H_{0}\left(  v_{n}^{+}\left(  \text{\textperiodcentered
}\right)  \right)  -4n\sum_{j=1}^{n-1}\delta_{\frac{j}{n}}\left(  v\left(
\text{\textperiodcentered}\right)  -v_{n}^{+}\left(  \text{\textperiodcentered
}\right)  \right)  ,\varphi\right\rangle _{\left(  H^{-1},H_{0}^{1}\right)
}\right\vert \\
& \leq C_{3}\left(  \left\Vert z\right\Vert _{C^{1}([0,1])}+\left\Vert
z\right\Vert _{C([0,1])}^{\frac{1}{2}}\right)  \left\Vert z\right\Vert
_{C([0,1])}\left\Vert \varphi\right\Vert _{H_{0}^{1}}.
\end{align*}
Thus, (\ref{AcotzGen}) follows.
\end{proof}

\bigskip

\begin{corollary}
There is $\rho\leq\delta$ such that if $v\in O_{\rho}(v_{n}^{+})$ and
$3/4<\alpha<1$, then
\begin{equation}
\left\Vert H_{0}\left(  v\left(  \text{\textperiodcentered}\right)  \right)
-H_{0}\left(  v_{n}^{+}\left(  \text{\textperiodcentered}\right)  \right)
-4n\sum_{j=1}^{n-1}\delta_{\frac{j}{n}}z\left(  \text{\textperiodcentered
}\right)  \right\Vert _{H^{-1}}\leq C\left\Vert z\right\Vert _{V_{H}^{2r}%
}^{\frac{1}{2}}\left\Vert z\right\Vert _{H_{0}^{1}},\label{Acotz1Gen}%
\end{equation}%
\begin{equation}
\left\Vert H_{0}\left(  v\left(  \text{\textperiodcentered}\right)  \right)
-H_{0}\left(  v_{n}^{+}\left(  \text{\textperiodcentered}\right)  \right)
-4n\sum_{j=1}^{n-1}\delta_{\frac{j}{n}}z\left(  \text{\textperiodcentered
}\right)  \right\Vert _{H^{-1}}\leq C\left\Vert z\right\Vert _{V_{H}^{2r}%
}^{\frac{1}{2}}\left\Vert z\right\Vert _{V^{2\alpha}},\label{Acotz2Gen}%
\end{equation}%
\begin{equation}
\left\Vert H_{0}\left(  v\left(  \text{\textperiodcentered}\right)  \right)
-H_{0}\left(  v_{n}^{+}\left(  \text{\textperiodcentered}\right)  \right)
-4n\sum_{j=1}^{n-1}\delta_{\frac{j}{n}}z\left(  \text{\textperiodcentered
}\right)  \right\Vert _{H^{-1}}\leq C\left\Vert z\right\Vert _{V_{H}^{2r}%
}^{\frac{3}{2}},\label{Acotz3Gen}%
\end{equation}
where $C>0$ is a constant.
\end{corollary}

\begin{proof}
It follows from (\ref{AcotzGen}) and the embeddings $V_{H}^{2r}\subset
C^{1}([0,1]),\ H_{0}^{1}(0,1)\subset C([0,1]),\ V^{2\alpha}\subset C([0,1])$
(see Lemma \ref{EmbeddingCLambda}).
\end{proof}

\bigskip

We will prove further that for any $n\geq2$ the function $g\left(  z\right)  $
given by
\[
g(z)=H_{0}(v_{n}^{\pm}+z)-H_{0}(v_{n}^{\pm})-4n\sum_{j=1}^{n-1}\delta
_{\frac{j}{n}}z\left(  \text{\textperiodcentered}\right)
\]
is Lipschitz in $V_{H}^{2r}$ in a suitable sense.

\begin{lemma}
\label{LipGen}There is $\rho\leq\delta$ such that if $v_{i}\in O_{\rho}%
(v_{n}^{\pm})$, then%
\begin{equation}
\left\Vert g\left(  z_{1}\left(  \text{\textperiodcentered}\right)  \right)
-g\left(  z_{2}\left(  \text{\textperiodcentered}\right)  \right)  \right\Vert
_{H^{-1}}\leq C(\left\Vert z_{1}\right\Vert _{V^{2r}}^{\frac{1}{2}}+\left\Vert
z_{2}\right\Vert _{V^{2r}}^{\frac{1}{2}})\left\Vert z_{1}-z_{2}\right\Vert
_{H_{0}^{1}},\label{gLipschitzH1Gen}%
\end{equation}%
\begin{equation}
\left\Vert g\left(  z_{1}\left(  \text{\textperiodcentered}\right)  \right)
-g\left(  z_{2}\left(  \text{\textperiodcentered}\right)  \right)  \right\Vert
_{H^{-1}}\leq C(\left\Vert z_{1}\right\Vert _{V^{2r}}^{\frac{1}{2}}+\left\Vert
z_{2}\right\Vert _{V^{2r}}^{\frac{1}{2}})\left\Vert z_{1}-z_{2}\right\Vert
_{V_{H}^{2r}},\label{gLipschitzVrGen}%
\end{equation}
for some $C>0$, where $z_{i}=v_{i}-v_{n}^{\pm}.$
\end{lemma}

\begin{proof}
Let us consider $v_{n}^{+}$. For $v_{1},v_{2}\in O_{\delta}\left(  v_{n}%
^{+}\right)  $ we have that%
\begin{align*}
& \left\langle g\left(  z_{1}\left(  \text{\textperiodcentered}\right)
\right)  -g\left(  z_{2}\left(  \text{\textperiodcentered}\right)  \right)
,\varphi\right\rangle \\
& =\left\langle H_{0}(u_{n}^{+}\left(  \text{\textperiodcentered}\right)
+z_{1}\left(  \text{\textperiodcentered}\right)  )-H_{0}(u_{n}^{+}\left(
\text{\textperiodcentered}\right)  )-4n\sum_{j=1}^{n-1}\delta_{\frac{j}{n}%
}z_{1}\left(  \text{\textperiodcentered}\right)  -(H_{0}(u_{n}^{+}\left(
\text{\textperiodcentered}\right)  +z_{2}\left(  \text{\textperiodcentered
}\right)  )-H_{0}(u_{n}^{+}\left(  \text{\textperiodcentered}\right)
)-4n\sum_{j=1}^{n-1}\delta_{\frac{j}{n}}z_{2}\left(  \text{\textperiodcentered
}\right)  ),\varphi\right\rangle \\
& =\left\langle H_{0}(v_{1}\left(  \text{\textperiodcentered}\right)
)-H_{0}(v_{2}\left(  \text{\textperiodcentered}\right)  )-4n\sum_{j=1}%
^{n-1}\delta_{\frac{j}{n}}(v_{1}\left(  \text{\textperiodcentered}\right)
-v_{2}\left(  \text{\textperiodcentered}\right)  ),\varphi\right\rangle
_{\left(  H^{-1},H_{0}^{1}\right)  }.
\end{align*}
Then%
\[
H_{0}\left(  v_{1}\left(  x\right)  \right)  -H_{0}\left(  v_{2}\left(
x\right)  \right)  =\left\{
\begin{array}
[c]{c}%
2(-1)^{j+1}\text{ for }x\in\left(  x_{v_{2}}^{j},x_{v_{1}}^{j}\right)  ,\text{
if }x_{v_{1}}^{j}>x_{v_{2}}^{j},\\
-2(-1)^{j+1}\text{ for }x\in\left(  x_{v_{1}}^{j},x_{v_{2}}^{j}\right)
,\text{ if }x_{v_{1}}^{j}<x_{v_{2}}^{j},\\
0\text{, otherwise,}%
\end{array}
\right.
\]
and for $\varphi\in H_{0}^{1}\left(  0,1\right)  ,$%
\[
\left\langle H_{0}\left(  v_{1}\left(  \text{\textperiodcentered}\right)
\right)  -H_{0}\left(  v_{2}\left(  \text{\textperiodcentered}\right)
\right)  ,\varphi\right\rangle _{\left(  H^{-1},H_{0}^{1}\right)  }=\sum
_{j=1}^{n-1}\int_{x_{v_{1}}^{j}}^{x_{v_{2}}^{j}}2(-1)^{j}\varphi\left(
x\right)  dx=\sum_{j=1}^{n-1}(-1)^{j}2\left(  x_{v_{2}}^{j}-x_{v_{1}}%
^{j}\right)  \varphi\left(  x_{j}^{\ast}\right)  ,
\]
for some $x_{_{j}}^{\ast}\in\lbrack\min\{x_{v_{1}}^{j},x_{v_{2}}^{j}%
\},\max\{x_{v_{1}}^{j},x_{v_{2}}^{j}\}]$. It follows that%
\begin{align*}
& \left\langle H_{0}(v_{1}\left(  \text{\textperiodcentered}\right)
)-H_{0}(v_{2}\left(  \text{\textperiodcentered}\right)  )-4n\sum_{j=1}%
^{n-1}\delta_{\frac{j}{n}}(v_{1}\left(  \text{\textperiodcentered}\right)
-v_{2}\left(  \text{\textperiodcentered}\right)  ),\varphi\right\rangle
_{\left(  H^{-1},H_{0}^{1}\right)  }\\
& =\sum_{j=1}^{n-1}(-1)^{j}2\left(  x_{v_{2}}^{j}-x_{v_{1}}^{j}\right)
\varphi\left(  x_{j}^{\ast}\right)  -4n\sum_{j=1}^{n-1}\left(  v_{1}\left(
\frac{j}{n}\right)  -v_{2}\left(  \frac{j}{n}\right)  \right)  \varphi\left(
\frac{j}{n}\right)  .
\end{align*}
Note that%
\begin{align*}
0  & =v_{1}(x_{v_{1}}^{j})-v_{2}(x_{v_{2}}^{j})=v_{1}(x_{v_{1}}^{j}%
)-v_{1}(x_{v_{2}}^{j})+v_{1}(x_{v_{2}}^{j})-v_{2}(x_{v_{2}}^{j})\\
& =v_{1}^{\prime}(\widetilde{x}_{j})(x_{v_{1}}^{j}-x_{v_{2}}^{j}%
)+v_{1}(x_{v_{2}}^{j})-v_{2}(x_{v_{2}}^{j}),
\end{align*}
for some $\widetilde{x}_{j}\in\lbrack\min\{x_{v_{1}}^{j},x_{v_{2}}^{j}%
\},\max\{x_{v_{1}}^{j},x_{v_{2}}^{j}\}]$. Hence,%
\begin{align*}
& \left\vert \left\langle H_{0}(v_{1}\left(  \text{\textperiodcentered
}\right)  )-H_{0}(v_{2}\left(  \text{\textperiodcentered}\right)
)-4n\sum_{j=1}^{n-1}\delta_{\frac{j}{n}}(v_{1}\left(
\text{\textperiodcentered}\right)  -v_{2}\left(  \text{\textperiodcentered
}\right)  ),\varphi\right\rangle _{\left(  H^{-1},H_{0}^{1}\right)
}\right\vert \\
& \leq\left\vert \sum_{j=1}^{n-1}\frac{2(-1)^{j}\varphi\left(  x_{j}^{\ast
}\right)  }{v_{1}^{\prime}(\widetilde{x}_{j})}(v_{1}(x_{v_{2}}^{j}%
)-v_{2}(x_{v_{2}}^{j}))-4n\sum_{j=1}^{n-1}\left(  v_{1}\left(  \frac{j}%
{n}\right)  -v_{2}\left(  \frac{j}{n}\right)  \right)  \varphi\left(  \frac
{j}{n}\right)  \right\vert \\
& \leq\sum_{j=1}^{n-1}2\left\vert \varphi\left(  x_{_{j}}^{\ast}\right)
\right\vert \left\vert \frac{(-1)^{j}}{v_{1}^{\prime}(\widetilde{x}_{_{j}}%
)}-2n\right\vert \left\vert v_{1}(x_{v_{2}}^{j})-v_{2}(x_{v_{2}}%
^{j})\right\vert +4n\sum_{j=1}^{n-1}\left\vert \varphi\left(  x_{j}^{\ast
}\right)  -\varphi\left(  \frac{j}{n}\right)  \right\vert \left\vert
v_{1}\left(  \frac{j}{n}\right)  -v_{2}\left(  \frac{j}{n}\right)  \right\vert
\\
& +4n\sum_{j=1}^{n-1}\left\vert \varphi\left(  x_{j}^{\ast}\right)
\right\vert \left\vert v_{1}(x_{v_{2}}^{j})-v_{2}(x_{v_{2}}^{j})-\left(
v_{1}\left(  \frac{j}{n}\right)  -v_{2}\left(  \frac{j}{n}\right)  \right)
\right\vert \\
& =B1+B2+B3.
\end{align*}
We estimate each term. For each $j$ we choose $k\left(  j\right)  \in\{1,2\}$
such that $x_{_{j}}^{\ast}\in\lbrack\min\{\frac{j}{n},x_{v_{k\left(  j\right)
}}^{j}\},\max\{\frac{j}{n},x_{v_{k\left(  j\right)  }}^{j}\}]$. By
\[
\left\vert \varphi\left(  x_{_{j}}^{\ast}\right)  -\varphi\left(  \frac{j}%
{n}\right)  \right\vert \leq\left\Vert \varphi\right\Vert _{H_{0}^{1}%
}\left\vert x_{_{j}}^{\ast}-\frac{j}{n}\right\vert ^{\frac{1}{2}}%
\]
and using
\[
-v_{k(j)}\left(  \frac{j}{n}\right)  =v_{k(j)}\left(  x_{v_{k\left(  j\right)
}}^{j}\right)  -v_{k(j)}\left(  \frac{j}{n}\right)  =v_{k(j)}^{\prime}\left(
w_{j}\right)  \left(  x_{v_{k\left(  j\right)  }}^{j}-\frac{j}{n}\right)  ,
\]
for some $w_{j}\in\lbrack\min\{\frac{j}{n},x_{v_{k\left(  j\right)  }}%
^{j}\},\max\{\frac{j}{n},x_{v_{k\left(  j\right)  }}^{j}\}]$, we get%
\begin{equation}
B_{2}\leq C_{1}\left\Vert \varphi\right\Vert _{H_{0}^{1}}\sum_{j=1}%
^{n-1}\left\vert v_{n}^{+}(\frac{j}{n})-v_{k(j)}(\frac{j}{n})\right\vert
^{\frac{1}{2}}\left\vert v_{1}\left(  \frac{j}{n}\right)  -v_{2}\left(
\frac{j}{n}\right)  \right\vert ,\label{B2Gen}%
\end{equation}
where we have used that there are $\gamma_{0}>0$ and $\rho\leq\delta$ such
that if $\left\Vert v_{2}-u_{n}^{+}\right\Vert _{V_{H}^{2r}}\leq\rho$, then
$\left\vert v_{k(j)}^{\prime}\left(  \widetilde{x}_{j}\right)  \right\vert
\geq\gamma_{0}$ for all $j$. Second, we have
\begin{align*}
B_{1}  & \leq C_{2}\sum_{j=1}^{n-1}\left\vert \varphi\left(  x_{j}^{\ast
}\right)  \right\vert \left(  \left\vert v_{1}^{\prime}\left(  \widetilde{x}%
_{j}\right)  -\left(  v_{n}^{+}\right)  ^{\prime}\left(  \widetilde{x}%
_{j}\right)  \right\vert +\left\vert \left(  v_{n}^{+}\right)  ^{\prime
}\left(  \widetilde{x}_{j}\right)  -\left(  v_{n}^{+}\right)  ^{\prime}\left(
\frac{j}{n}\right)  \right\vert \right)  \left\vert v_{1}(x_{v_{2}}^{j}%
)-v_{2}(x_{v_{2}}^{j})\right\vert \\
& \leq C_{3}\sum_{j=1}^{n-1}\left\vert \varphi\left(  x_{j}^{\ast}\right)
\right\vert \left(  \left\vert v_{1}^{\prime}\left(  \widetilde{x}_{j}\right)
-\left(  v_{n}^{+}\right)  ^{\prime}\left(  \widetilde{x}_{j}\right)
\right\vert +\left\Vert v_{n}^{+}\right\Vert _{H^{2}}\left\vert \widetilde{x}%
_{j}-\frac{j}{n}\right\vert ^{\frac{1}{2}}\right)  \left\vert v_{1}(x_{v_{2}%
}^{j})-v_{2}(x_{v_{2}}^{j})\right\vert .
\end{align*}
For each $j$ we choose $r\left(  j\right)  \in\{1,2\}$ such that
$\widetilde{x}_{j}\in\lbrack\min\{\frac{j}{n},x_{v_{r\left(  j\right)  }}%
^{j}\},\max\{\frac{j}{n},x_{v_{r\left(  j\right)  }}^{j}\}]$. Then
\[
-v_{r(j)}\left(  \frac{j}{n}\right)  =v_{r(j)}\left(  x_{v_{r(j)}}^{j}\right)
-v_{r(j)}\left(  \frac{j}{n}\right)  =v_{r(j)}^{\prime}\left(  b_{j}\right)
\left(  x_{v_{r(j)}}^{j}-\frac{j}{n}\right)  ,
\]
for some $b_{j}\in\lbrack\min\{x_{v_{r\left(  j\right)  }}^{j},\frac{j}%
{n}\},\max\{x_{v_{r\left(  j\right)  }}^{j},\frac{j}{n}\}]$, so%
\begin{equation}
B_{1}\leq C_{4}\sum_{j=1}^{n-1}\left\vert \varphi\left(  x_{j}^{\ast}\right)
\right\vert (\left\vert v_{1}^{\prime}(\widetilde{x}_{j})-\left(  v_{n}%
^{+}\right)  ^{\prime}(\widetilde{x}_{j})\right\vert +\left\Vert v_{n}%
^{+}\right\Vert _{H^{2}}\left\vert v_{n}^{+}(\frac{j}{n})-v_{r(j)}(\frac{j}%
{n})\right\vert ^{\frac{1}{2}})\left\vert v_{1}(x_{v_{2}}^{j})-v_{2}(x_{v_{2}%
}^{j})\right\vert .\label{B1Gen}%
\end{equation}
Further, by%
\begin{align*}
& \left\vert v_{1}(x_{v_{2}}^{j})-v_{2}(x_{v_{2}}^{j})-\left(  v_{1}\left(
\frac{j}{n}\right)  -v_{2}\left(  \frac{j}{n}\right)  \right)  \right\vert \\
& =\left\vert \int_{\frac{j}{n}}^{x_{v_{2}}}v_{1}^{\prime}(s)ds-\int_{\frac
{j}{n}}^{x_{v_{2}}}v_{2}^{\prime}(s)ds\right\vert \\
& \leq\left\Vert v_{1}-v_{2}\right\Vert _{H_{0}^{1}}\left\vert x_{v_{2}}%
-\frac{j}{n}\right\vert ^{\frac{1}{2}}%
\end{align*}
and using (\ref{Eqvj}) and $\left\vert v_{2}^{\prime}\left(  \widetilde{x}%
_{v_{2},j}\right)  \right\vert \geq\gamma_{0}$ we obtain%
\begin{equation}
B_{3}\leq C_{5}\left\Vert \varphi\right\Vert _{H_{0}^{1}}\left\vert v_{n}%
^{+}(\frac{j}{n})-v_{2}(\frac{j}{n})\right\vert ^{\frac{1}{2}}\left\Vert
v_{1}-v_{2}\right\Vert _{H_{0}^{1}}.\label{B3Gen}%
\end{equation}
Hence, (\ref{B2Gen})-(\ref{B3Gen}) imply that%
\begin{align*}
& \left\vert \left\langle H_{0}(v_{1}\left(  \text{\textperiodcentered
}\right)  )-H_{0}(v_{2}\left(  \text{\textperiodcentered}\right)
)-4n\sum_{j=1}^{n-1}\delta_{\frac{j}{n}}(v_{1}\left(
\text{\textperiodcentered}\right)  -v_{2}\left(  \text{\textperiodcentered
}\right)  ),\varphi\right\rangle _{\left(  H^{-1},H_{0}^{1}\right)
}\right\vert \\
& \leq C_{6}\left\Vert \varphi\right\Vert _{H_{0}^{1}}(\left\Vert v_{1}%
-v_{n}^{+}\right\Vert _{C^{1}}^{\frac{1}{2}}+\left\Vert v_{2}-v_{n}%
^{+}\right\Vert _{C^{1}}^{\frac{1}{2}})\left\Vert v_{1}-v_{2}\right\Vert
_{H_{0}^{1}}\\
& \leq C_{7}\left\Vert \varphi\right\Vert _{H_{0}^{1}}(\left\Vert v_{1}%
-v_{n}^{+}\right\Vert _{V^{2r}}^{\frac{1}{2}}+\left\Vert v_{2}-v_{n}%
^{+}\right\Vert _{V^{2r}}^{\frac{1}{2}})\left\Vert v_{1}-v_{2}\right\Vert
_{H_{0}^{1}}\\
& \leq C_{8}\left\Vert \varphi\right\Vert _{H_{0}^{1}}(\left\Vert v_{1}%
-v_{n}^{+}\right\Vert _{V^{2r}}^{\frac{1}{2}}+\left\Vert v_{2}-v_{n}%
^{+}\right\Vert _{V^{2r}}^{\frac{1}{2}})\left\Vert v_{1}-v_{2}\right\Vert
_{V_{H}^{2r}},
\end{align*}
so (\ref{gLipschitzH1Gen})-(\ref{gLipschitzVrGen}) follow.
\end{proof}

\bigskip

Let us consider a solution $u\left(  \text{\textperiodcentered}\right)  $ to
problem (\ref{Eq}) such that $u\left(  t\right)  \in O_{\delta}(v_{n}^{\pm})$,
$n\geq2$, for $t\in\lbrack0,t_{0}]$, being $\rho>0$ the constant from Lemma
\ref{AcotzLemma2}. We know that it is unique on $[0,t_{0}]$. The difference
$z=u-v_{n}^{\pm}$ satisfies%
\begin{equation}
\frac{dz}{dt}+Az\left(  t\right)  =H_{0}\left(  u\left(  t\right)  \right)
-H_{0}\left(  v_{n}^{\pm}\right)  .\label{z}%
\end{equation}
Hence,%
\[
\frac{dz}{dt}+Lz\left(  t\right)  =H_{0}\left(  u\left(  t\right)  \right)
-H_{0}\left(  v_{n}^{\pm}\right)  -4n\sum_{k=1}^{n-1}\delta_{\frac{k}{n}%
}z\left(  \text{\textperiodcentered}\right)  =g\left(  z\left(  t\right)
\right)  .
\]
Putting $z\left(  t\right)  =E_{1}z\left(  t\right)  +E_{2}(t)=z_{1}%
(t)+z_{2}(t)$, we have%
\begin{align}
\frac{dz_{1}}{dt}+L_{1}z_{1}\left(  t\right)   & =E_{1}g\left(  z\left(
t\right)  \right)  ,\label{z1}\\
\frac{dz_{2}}{dt}+L_{2}z_{2}\left(  t\right)   & =E_{2}g\left(  z\left(
t\right)  \right)  ,\label{z2}%
\end{align}
and $z_{i}(0)=E_{i}z\left(  0\right)  $, where the equalities are satisfied
for a.a. $t\in\left(  0,t_{0}\right)  $ in the spaces $V_{1}^{\prime}$ and
$V_{2}^{\prime}$, respectively. Since $z\left(  \text{\textperiodcentered
}\right)  \in C([0,\infty),H_{0}^{1}\left(  0,1\right)  )$, it is easy to see
that $g\left(  z\left(  \text{\textperiodcentered}\right)  \right)  \in
L_{loc}^{\infty}(0,\infty;H^{-1}(0,1))$. These functions satisfy the formula
of variation of constants:%
\begin{align}
z_{1}(t)  & =e^{-L_{1}(t-l)}z_{1}\left(  l\right)  +\int_{l}^{t}%
e^{-L_{1}\left(  t-s\right)  }E_{1}g\left(  z(s)\right)  ds,\label{VC}\\
z_{2}(t)  & =e^{-L_{2}(t-l)}z_{2}\left(  l\right)  +\int_{l}^{t}%
e^{-L_{2}\left(  t-s\right)  }E_{2}g\left(  z(s)\right)  ds,\ 0\leq l\leq
t\leq t_{0}.\nonumber
\end{align}
Indeed, since $z\left(  \text{\textperiodcentered}\right)  $ is a strong
solution to (\ref{z}), it is easy to see that the functions $z_{1}\left(
\text{\textperiodcentered}\right)  ,\ z_{2}\left(  \text{\textperiodcentered
}\right)  $ are, respectively, strong solutions to problems (\ref{z1}%
)-(\ref{z2}) in the spaces $V_{1}^{\prime}$, $V_{2}^{\prime}$ endowed with the
induced topology from $H^{-1}(0,1)$. As $E_{i}g\left(  z\left(
\text{\textperiodcentered}\right)  \right)  \in L_{loc}^{\infty}%
(0,\infty;V_{i}^{\prime})$, Lemmas \ref{MildLemma}, \ref{DisL1L2} imply
(\ref{VC}).

\begin{theorem}
\label{Unstable2rGen}The fixed points $v_{n}^{+}$ ($v_{n}^{-}$), $n\geq2$, are
unstable in $V_{H}^{2r}.$
\end{theorem}

\begin{proof}
By contradiction assume that $v_{n}^{+}$ ($v_{n}^{-}$) is stable. Then, for
any constant $\delta_{0}<\rho$ (where $\rho$ is defined in Lemma
\ref{AcotzLemma2}) there exists $\gamma>0$ such that if $\left\Vert z\left(
0\right)  \right\Vert _{V_{H}^{2r}}\leq\gamma$, then the unique solution
$u\left(  \text{\textperiodcentered}\right)  $ to problem (\ref{Eq}) satisfies
$\left\Vert z\left(  t\right)  \right\Vert _{V_{H}^{2r}}\leq\delta_{0}$ for
any $t\geq0$.

For $t>0$ we denote $\widetilde{z}_{1}(s)=z_{1}(t-s)$,$\ \widetilde{z}%
(s)=z(t-s),\ s\in\lbrack0,t]$. Then by (\ref{VC}),%
\begin{align*}
\widetilde{z}_{1}(0)  & =z_{1}(t)=e^{-L_{1}s}z_{1}\left(  t-s\right)
+\int_{t-s}^{t}e^{-L_{1}\left(  t-l\right)  }E_{1}g\left(  z(l)\right)  dl\\
& =e^{-L_{1}s}\widetilde{z}_{1}\left(  s\right)  +\int_{0}^{s}e^{-L_{1}r}%
E_{1}g\left(  \widetilde{z}(r)\right)  dr,
\end{align*}
so%
\[
\widetilde{z}_{1}\left(  s\right)  =e^{L_{1}s}\widetilde{z}_{1}(0)-\int%
_{0}^{s}e^{L_{1}(s-r)}E_{1}g\left(  \widetilde{z}(r)\right)  dr
\]
and applying (\ref{IneqL2C})-(\ref{IneqL2D}) with $\alpha\in\left(
3/4,1\right)  $ and (\ref{Acotz2Gen}) we get
\begin{align*}
\left\Vert \widetilde{z}_{1}\left(  s\right)  \right\Vert _{V^{2\alpha}}  &
\leq C_{1}e^{\lambda_{m}s}\left\Vert \widetilde{z}_{1}\left(  0\right)
\right\Vert _{V^{2\alpha}}+C_{2}\int_{0}^{s}e^{\lambda_{m}(s-r)}\left\Vert
E_{1}g\left(  \widetilde{z}(r)\right)  \right\Vert _{H^{-1}}dr\\
& \leq C_{1}e^{\lambda_{m}s}\left\Vert \widetilde{z}_{1}\left(  0\right)
\right\Vert _{V^{2\alpha}}+C_{3}\left\Vert E_{1}\right\Vert \int_{0}%
^{s}e^{\lambda_{m}(s-r)}\left\Vert \widetilde{z}(r)\right\Vert _{V_{H}^{2r}%
}^{\frac{1}{2}}\left\Vert \widetilde{z}(r)\right\Vert _{V^{2\alpha}}dr\\
& \leq C_{1}e^{\lambda_{m}s}\left\Vert \widetilde{z}_{1}\left(  0\right)
\right\Vert _{V^{2\alpha}}+C_{4}\delta_{0}^{1/2}\int_{0}^{s}e^{\lambda
_{m}(s-r)}\left\Vert \widetilde{z}(r)\right\Vert _{V^{2\alpha}}dr.
\end{align*}

If $\left\Vert z_{2}(0)\right\Vert _{V^{2\alpha}}<\left\Vert z_{1}%
(0)\right\Vert _{V^{2\alpha}}$, then there is $t_{1}>0$ such that $\left\Vert
z_{2}(t)\right\Vert _{V^{2\alpha}}\leq\left\Vert z_{1}(t)\right\Vert
_{V^{2\alpha}}$ if $t\in\lbrack0,t_{1}]$. Putting $y\left(  s\right)
=e^{-\lambda_{m}s}\left\Vert \widetilde{z}_{1}\left(  s\right)  \right\Vert
_{V^{2\alpha}}$ and using $\left\Vert \widetilde{z}(s)\right\Vert
_{V^{2\alpha}}\leq2\left\Vert \widetilde{z}_{1}(s)\right\Vert _{V^{2\alpha}}$
we find that%
\[
y\left(  s\right)  \leq C_{1}\left\Vert \widetilde{z}_{1}\left(  0\right)
\right\Vert _{V^{2\alpha}}+2C_{4}\delta_{0}^{1/2}\int_{0}^{s}y\left(
r\right)  dr,\text{ for }s\in\lbrack0,t],
\]
and Gronwall's lemma implies%
\[
y\left(  s\right)  =e^{-\lambda_{m}s}\left\Vert \widetilde{z}_{1}\left(
s\right)  \right\Vert _{V^{2\alpha}}\leq C_{1}\left\Vert \widetilde{z}%
_{1}\left(  0\right)  \right\Vert _{V^{2\alpha}}e^{2C_{4}\delta_{0}^{1/2}%
s}\text{ }\forall s\in\lbrack0,t].
\]
Thus,%
\begin{equation}
\sup_{s\in\lbrack0,t]}\left\Vert z_{1}\left(  t-s\right)  \right\Vert
_{V^{2\alpha}}e^{\left(  -\lambda_{m}-2C_{4}\delta_{0}^{1/2}\right)  s}\leq
C_{1}\left\Vert z_{1}\left(  t\right)  \right\Vert _{V^{2\alpha}%
}.\label{Ineqz1}%
\end{equation}
We choose $\delta_{0}$ such that $-\lambda_{m}-2C_{4}\delta_{0}^{1/2}%
\geq-\lambda_{m}/2.$

On the other hand, using (\ref{VC}), (\ref{IneqL2A})-(\ref{IneqL2B}),
(\ref{Acotz2Gen}) and $\left\Vert z(s)\right\Vert _{V^{2\alpha}}%
\leq2\left\Vert z_{1}(s)\right\Vert _{V^{2\alpha}}$ we obtain%
\begin{align*}
\left\Vert z_{2}(t)\right\Vert _{V^{2\alpha}}  & \leq C_{5}e^{-\lambda_{m+1}%
t}\left\Vert z_{2}\left(  0\right)  \right\Vert _{V^{2\alpha}}+C_{6}\int%
_{0}^{t}e^{-\lambda_{m+1}\left(  t-s\right)  }\left(  t-s\right)  ^{-\alpha
}\left\Vert E_{2}g\left(  z(s)\right)  \right\Vert _{H^{-1}}ds\\
& \leq C_{5}e^{-\lambda_{m+1}t}\left\Vert z_{2}\left(  0\right)  \right\Vert
_{V^{2\alpha}}+2C_{7}\left\Vert E_{2}\right\Vert \int_{0}^{t}e^{-\lambda
_{m+1}\left(  t-s\right)  }\left(  t-s\right)  ^{-\alpha}\left\Vert
z(s)\right\Vert _{V_{H}^{2r}}^{\frac{1}{2}}\left\Vert z_{1}(r)\right\Vert
_{V^{2\alpha}}ds\\
& \leq C_{5}e^{-\lambda_{m+1}t}\left\Vert z_{2}\left(  0\right)  \right\Vert
_{V^{2\alpha}}+C_{8}\delta_{0}^{1/2}\sup_{s\in\lbrack0,t]}\left\Vert
z_{1}\left(  t-s\right)  \right\Vert _{V^{2\alpha}}\int_{0}^{\infty
}e^{-\lambda_{m+1}r}r^{-\alpha}dr.
\end{align*}
Hence,%
\begin{equation}
C_{1}\left\Vert z_{2}(t)\right\Vert _{V^{2\alpha}}\leq C_{9}\left\Vert
z_{2}\left(  0\right)  \right\Vert _{V^{2\alpha}}+C_{10}\delta_{0}^{1/2}%
\sup_{s\in\lbrack0,t]}\left\Vert z_{1}\left(  t-s\right)  \right\Vert
_{V^{2\alpha}}e^{\left(  -\lambda_{m}-2C_{4}\delta_{0}^{1/2}\right)
s},\label{Ineqz2}%
\end{equation}
where we have used that $-\lambda_{m}-2C_{4}\delta_{0}^{1/2}>0$, so
$e^{\left(  -\lambda_{m}-2C_{4}\delta_{0}^{1/2}\right)  s}>1$. By
(\ref{Ineqz1}), (\ref{Ineqz2}) we have%
\[
C_{1}\left(  \left\Vert z_{1}\left(  t\right)  \right\Vert _{V^{2\alpha}%
}-\left\Vert z_{2}\left(  t\right)  \right\Vert _{V^{2\alpha}}\right)
\geq\sup_{s\in\lbrack0,t]}\left\Vert z_{1}\left(  t-s\right)  \right\Vert
_{V^{2\alpha}}e^{\left(  -\lambda_{m}-2C_{4}\delta_{0}^{1/2}\right)  s}\left(
1-C_{10}\delta_{0}^{1/2}\right)  -C_{9}\left\Vert z_{2}\left(  0\right)
\right\Vert _{V^{2\alpha}}.
\]
We choose $\delta_{0}$ such that $1-C_{10}\delta_{0}^{1/2}\geq1/2,\ -\lambda
_{m}-2C_{4}\delta_{0}^{1/2}\geq-\lambda_{m}/2.$ Then,%
\[
C_{1}\left(  \left\Vert z_{1}\left(  t\right)  \right\Vert _{V^{2\alpha}%
}-\left\Vert z_{2}\left(  t\right)  \right\Vert _{V^{2\alpha}}\right)
\geq\frac{1}{2}\left\Vert z_{1}\left(  0\right)  \right\Vert _{V^{2\alpha}%
}e^{-\frac{\lambda_{m}}{2}t}-C_{9}\left\Vert z_{2}\left(  0\right)
\right\Vert _{V^{2\alpha}}.
\]
Then if $C_{9}\left\Vert z_{2}\left(  0\right)  \right\Vert _{V^{2\alpha}}%
\leq(1/2)\left\Vert z_{1}\left(  0\right)  \right\Vert _{V^{2\alpha}}$, we
deduce that $\left\Vert z_{2}(t)\right\Vert _{V^{2\alpha}}\leq\left\Vert
z_{1}(t)\right\Vert _{V^{2\alpha}}$ for any $t\geq0$. It follows from
(\ref{Ineqz1}) that%
\begin{equation}
\left\Vert z_{1}\left(  t\right)  \right\Vert _{V^{2\alpha}}\geq\frac{1}%
{C_{1}}\left\Vert z_{1}\left(  0\right)  \right\Vert _{V^{2\alpha}}%
e^{-\frac{\lambda_{m}}{2}t}\text{ }\forall t\geq0\text{. }\label{Ineqz1B}%
\end{equation}
The embedding $V_{H}^{2r}\subset V^{2\alpha}$ implies that%
\[
\left\Vert z_{1}\left(  t\right)  \right\Vert _{V^{2\alpha}}\leq
K_{1}\left\Vert z\left(  t\right)  \right\Vert _{V^{2\alpha}}\leq
K_{2}\left\Vert z\left(  t\right)  \right\Vert _{V_{H}^{2r}}\leq K_{2}%
\delta_{0}\text{, }\forall t\geq0.
\]
But if $\delta_{0}<\delta_{1}$, from (\ref{Ineqz1B}) there is $\overline{t}>0$
such that
\[
\left\Vert z_{1}\left(  \overline{t}\right)  \right\Vert _{V^{2\alpha}}%
=K_{2}\delta_{1}>K_{2}\delta_{0},
\]
which is a contradiction.
\end{proof}

\bigskip

We recall that the function $\phi:\mathbb{R}\rightarrow L^{2}\left(
0,1\right)  $ is a complete trajectory of problem (\ref{Eq}) if $\phi\left(
\text{\textperiodcentered}+s\right)  $ is a strong solution for any
$s\in\mathbb{R}$.

It is well known \cite{ArrRBVal} that the multivalued semflow generated by the
solutions of problem (\ref{Eq}) possesses a global attractor consisting of all
bounded complete trajectories. Moreover, any bounded complete trajectory
converges backwards and forward to a fixed point in $L^{2}(0,1)$ (and in fact
in $V_{H}^{2r}$ as the global attractor is compact in $V_{H}^{2r}$ \cite[Lemma
14]{Valero2021}). The attractor is described by the fixed points and the
heteroclinic connections between them. As mentioned in the introduction, a
partial answer about which connections exist was given in \cite{ArrRBVal}. In
particular, homoclinic connections are forbidden, that is, there is no
complete bounded trajectories $\phi\left(  \text{\textperiodcentered}\right)
$ satisfying $\phi\left(  t\right)  \overset{t\rightarrow+\infty}{\rightarrow
}z,\ \phi\left(  t\right)  \overset{t\rightarrow-\infty}{\rightarrow}z$, where
$z$ is a fixed point. This implies that if a bounded complete trajectory
$\phi\left(  \text{\textperiodcentered}\right)  $ is contained in a closed
region $U$ where there is only one fixed point $z$, then $\phi\left(
t\right)  \equiv z$.

Also, it is known \cite[Lemma 31]{CostaValero} that if $u^{n}\left(
\text{\textperiodcentered}\right)  $ is a sequence of strong solutions such
that $u^{n}\left(  0\right)  \rightarrow u_{0}$, then up to a subsequence
$u^{n}\rightarrow u$ in $C([0,T],L^{2}\left(  0,1\right)  )$, for all $T>0$,
where $u\left(  \text{\textperiodcentered}\right)  $ is a strong solution
satisfying $u\left(  0\right)  =u_{0}.$

\begin{theorem}
The fixed points $v_{n}^{+}$ ($v_{n}^{-}$), $n\geq2$, are unstable in
$L^{2}(0,1).$
\end{theorem}

\begin{proof}
Let us consider $v_{n}^{+}$. By contradiction let $v_{n}^{+}$ be stable in
$L^{2}\left(  0,1\right)  $. Then for any $\varepsilon>0$ there is $\gamma>0$
such that if $\left\Vert u\left(  0\right)  -v_{n}^{+}\right\Vert _{L^{2}%
}<\gamma,$ then any solution $u\left(  \text{\textperiodcentered}\right)  $
starting at $u\left(  0\right)  $ satisfies $\left\Vert u\left(  t\right)
-v_{n}^{+}\right\Vert _{L^{2}}<\varepsilon$ for any $t\geq0.$ By Theorem
\ref{Unstable2rGen} there exist a constant $\beta>0,$ a sequence of solutions
$\{u_{n}($\textperiodcentered$)\}$ and a sequence of times $\{t_{n}\}$ such
that $u_{n}(0)\rightarrow v_{n}^{+}$ in $V_{H}^{2r}$, so $u_{n}(0)\rightarrow
v_{n}^{+}$ in $L^{2}(0,1)$ as well, and $\left\Vert u_{n}(t_{n})-v_{n}%
^{+}\right\Vert _{V_{H}^{2r}}>\beta$.

We will show that this is not possible.

First, let us show that the sequence $\{u_{n}(t_{n})\}$ is relatively compact
in $V_{H}^{2r}$ if $t_{n}\geq b$ for some $b>0$. By the variation of constants
formula (see Lemma 4 in \cite{CLV20}) we have%
\[
u_{n}(t)=e^{-A_{H}t}u_{n}(0)+\int_{0}^{t}e^{-A(t-s)}f_{n}(s)ds,
\]
where $f_{n}\in L_{loc}^{2}(0,\infty;L^{2}(0,1))$ are such that $f\left(
s,x\right)  \in H_{0}(u_{n}(s,x))$ for a.a. $x\in\left(  0,1\right)  ,\ s>0$.
It is easy to see that $\left\Vert f_{n}(t)\right\Vert _{L^{\infty}%
(0,\infty;L^{2}(0,1))}\leq C_{1}$. Thus, by the estimate of the norm of
$e^{-A_{H}t}$ in $V^{2\eta}$, $\eta\in\lbrack0,1)$ (see Theorem 37.5 in
\cite{SellYou}), there exist $a>0,\ M_{\eta}>0$ such that%
\begin{align}
\left\Vert u_{n}(t_{n})\right\Vert _{V_{H}^{2\eta}}  & \leq\left\Vert
A_{H}^{\eta}e^{-A_{H}t_{n}}u_{n}(0)\right\Vert _{L^{2}}+\int_{0}^{t_{n}%
}\left\Vert A_{H}^{\eta}e^{-A_{H}t_{n}}f_{n}(s)\right\Vert _{L^{2}%
}ds\nonumber\\
& \leq M_{\eta}e^{-at_{n}}t_{n}^{-\eta}\left\Vert u_{n}(0)\right\Vert _{L^{2}%
}+C_{1}M_{\eta}\int_{0}^{t_{n}}e^{-a\left(  t_{n}-s\right)  }\left(
t-s\right)  ^{-\eta}ds\leq C_{2},\label{IneqEta}%
\end{align}
so $\{u_{n}(t_{n})\}$ is bounded in $V_{H}^{2\eta}$. The compact embedding
$V_{H}^{2\eta}\subset V_{H}^{2r}$ for $r<\eta$ implies the result.

First, let $\{t_{n}\}$ be bounded, so we can assume that $t_{n}\rightarrow
t_{0}$. If $t_{0}>0$, then $t_{n}\geq t_{0}/2>0$. By the continuity property
described above, up to a subsequence $u_{n}(t_{n})\rightarrow u(t_{0})$ in
$L^{2}\left(  0,1\right)  $, where $u\left(  \text{\textperiodcentered
}\right)  $ is a solution with $u\left(  0\right)  =v_{n}^{+}$. As
$\{u_{n}(t_{n})\}$ is relatively compact in $V_{H}^{2r}$, $u_{n}%
(t_{n})\rightarrow u(t_{0})$ in $V_{H}^{2r}$ as well. Since $u(t)\equiv
v_{n}^{+}$ is the unique solution such that $u\left(  0\right)  =v_{n}^{+}$,
we deduce that $u(t_{0})=v_{n}^{+}$, which is a contradiction with $\left\Vert
u_{n}(t_{n})-v_{n}^{+}\right\Vert _{V_{H}^{2r}}>\beta$. If $t_{0}=0$, then
using again the variation of constants formula we get%
\begin{align*}
\left\Vert u_{n}(t_{n})-v_{n}^{+}\right\Vert _{V_{H}^{2r}}  & \leq\left\Vert
A_{H}^{r}e^{-A_{H}t_{n}}\left(  u_{n}(0)-v_{n}^{+}\right)  \right\Vert
_{L^{2}}+\int_{0}^{t_{n}}\left\Vert A_{H}^{r}e^{-A_{H}t_{n}}(f_{n}%
(s)-g_{n}(s)\right\Vert _{L^{2}}ds\\
& \leq C_{2}\left\Vert u_{n}(0)-v_{n}^{+}\right\Vert _{V_{H}^{2r}}+C_{3}%
\int_{0}^{t_{n}}e^{-a\left(  t_{n}-s\right)  }\left(  t-s\right)
^{-r}ds\rightarrow0\text{ as }n\rightarrow\infty,
\end{align*}
which contradicts $\left\Vert u_{n}(t_{n})-v_{n}^{+}\right\Vert _{V_{H}^{2r}%
}>\beta$.

Second, let $t_{n}\rightarrow\infty$. As $\{u_{n}(t_{n})\}$ is relatively
compact in $V_{H}^{2r}$, passing to a subsequence $u_{n}(t_{n})\rightarrow
y^{0}$ in $V_{H}^{2r}$. Denote $v_{n}^{0}\left(  \text{\textperiodcentered
}\right)  =u_{n}\left(  \text{\textperiodcentered}+t_{n}\right)  $. Then up to
a subsequence we get that $v_{n}^{0}\left(  t\right)  \rightarrow v^{0}(t)$ in
$L^{2}(0,1)$, for any $t\geq0$, where $v^{0}($\textperiodcentered$)$ is a
solution such that $v^{0}\left(  0\right)  =y^{0}$. It follows from the
stability assumption that $\left\Vert v^{0}(t)-v_{n}^{+}\right\Vert _{L^{2}%
}\leq\varepsilon$ for any $t\geq0$. Further, $v_{n}^{-1}\left(
\text{\textperiodcentered}\right)  =u_{n}\left(  \text{\textperiodcentered
}+t_{n}-1\right)  $. The sequence $\{u_{n}(t_{n}-1)\}$ is relatively compact
in $V_{H}^{2r}$, so $u_{n}(t_{n}-1)\rightarrow y^{-1}$ and thus we have
$v_{n}^{-1}\left(  t\right)  \rightarrow v^{-1}(t)$, for any $t\geq0$, where
$v^{-1}($\textperiodcentered$)$ is a solution such that $v^{-1}\left(
0\right)  =y^{-1}$ and $\left\Vert v^{-1}(t)-v_{n}^{+}\right\Vert _{L^{2}}%
\leq\varepsilon$ for any $t\geq0$. Moreover, $v^{-1}(t+1)=v^{0}\left(
t\right)  $ for $t\geq0$. In this way, we define a sequence of solutions
$v^{-k}\left(  \text{\textperiodcentered}\right)  $, $k\geq0$, such that
$v^{-k}(t+1)=v^{-k+1}(t)$ and $\left\Vert v^{-k}(t)-v_{n}^{+}\right\Vert
_{L^{2}}\leq\varepsilon$ for $t\geq0$. Define $\phi:\mathbb{R}\rightarrow
L^{2}(0,1)$ by $\phi(t)=v^{-k}(t+k)$ if $t\geq-k$. Then $\phi\left(
\text{\textperiodcentered}\right)  $ is a bounded complete trajectory. We
choose the constant $\varepsilon>0$ in such a way that in the ball
$O_{2\varepsilon}(v_{n}^{+})$ of the space $L^{2}(0,1)$ there is no other
fixed point. Hence, $\phi\equiv v_{n}^{+}$, which is a contradiction with
\[
\left\Vert \phi\left(  0\right)  -v_{n}^{+}\right\Vert _{V_{H}^{2r}%
}=\left\Vert y^{0}-v_{n}^{+}\right\Vert _{V_{H}^{2r}}=\lim_{n\rightarrow
\infty}\left\Vert u_{n}(t_{n})-v_{n}^{+}\right\Vert _{V_{H}^{2r}}\geq\beta.
\]

\end{proof}

\begin{remark}
In \cite{ArrRBVal} the instability of the fixed points $v_{n}^{\pm}$, $n\geq
2$, in the space $L^{2}\left(  0,1\right)  $ was proved using a Lyapunov
function. Here, using a linearization technique, we have proved their
instability in the more regular space $V_{H}^{2r}$, providing also an
alternative proof of instability in the space $L^{2}(0,1).$
\end{remark}

We finish this section by analyzing the stability of the fixed points
$v_{1}^{\pm}.$ In \cite{ArrRBVal} it was establish that $v_{1}^{\pm}$ are
asymptotically stable in the spaces $L^{2}(0,1)$ and $H_{0}^{1}(0,1)$. In the
following theorem we prove thei stability in the more regular space
$V_{H}^{2r}.$

\begin{theorem}
The fixed points $v_{1}^{+}$,$v_{1}^{-}$ are asymptotically stable in
$V_{H}^{2r}.$
\end{theorem}

\begin{proof}
We observe that
\[
\left\Vert e^{-A_{H}t}u_{0}\right\Vert _{V_{H}^{2r}}\leq Me^{-\mu_{1}%
t}\left\Vert u_{0}\right\Vert _{V_{H}^{2r}}\text{, for }t\geq0,
\]
where $M\geq1$ and $\mu_{1}>0$ denotes the first eigenvalue of the operator
$A_{H}$.

We choose an initial condition $u_{0}\in O_{\delta_{0}}(v_{1}^{\pm})$, where
$\delta_{0}<\delta$ satisfies $\delta_{0}\leq\frac{\delta}{2M}$. Let $z\left(
t\right)  =u\left(  t\right)  -v_{1}^{\pm}$, where $u\left(
\text{\textperiodcentered}\right)  $ is the solution corresponding to the
initial condition $u_{0}$, which is unique as long a $u\left(  t\right)  \in
O_{\delta}(v_{1}^{\pm})$. Also, we can choose $\delta>0$ such that $v\left(
x\right)  >0$, if $x\in\left(  0,1\right)  $, for any $v\in O_{\delta}%
(v_{1}^{\pm})$. Then, we have that%
\[
H_{0}\left(  v\left(  \text{\textperiodcentered}\right)  \right)
-H_{0}\left(  v_{1}^{\pm}\right)  \equiv0,
\]
on some maximal interval $[0,t_{1})$ we obtain that%
\[
z\left(  t\right)  =e^{-A_{H}t}z_{0},
\]%
\[
\left\Vert e^{-A_{H}t}z_{0}\right\Vert _{V_{H}^{2r}}\leq Me^{-\mu_{1}%
t}\left\Vert z_{0}\right\Vert _{V_{H}^{2r}}\leq\frac{\delta}{2}.
\]
The last inequalities imply that $t_{1}=+\infty$ and that $v_{1}^{+}$%
,$v_{1}^{-}$ are asymptotically stable in $V_{H}^{2r}.$
\end{proof}

\section{Saddle-point property}

In this section, we will prove that the stable and unstable sets of the fixed
points $v_{n}^{\pm}$, $n\geq2,$ are tangent to the linear subspaces generated,
respectively, by the negative and positive eigenvalues of the linearized
equation. This implies that they satisfy the saddle-point property, so they
are hyperbolic.

Let $0<\delta_{0}\leq\delta$, where $\delta>0$ is defined before, be such that
the $2\delta_{0}$-neighborhood $O_{2\delta_{0}}(v_{n}^{\pm})$ of $v_{n}^{\pm}$
in the space $V_{H}^{2r}$ does not contain any other fixed point. We recall
that for any $u_{0}$ in the neighborhood $O_{\delta}\left(  v_{n}^{\pm
}\right)  $ in $V_{H}^{2r}$, $3/4<r<1,$ there is a unique solution $u\left(
\text{\textperiodcentered}\right)  $. For a solution $u\left(
\text{\textperiodcentered}\right)  $ denote as before $z\left(
\text{\textperiodcentered}\right)  =u\left(  \text{\textperiodcentered
}\right)  -v_{n}^{\pm}$. We define the $\delta_{0}$-stable and unstable sets
of $v_{n}^{\pm}$ in $V_{H}^{2r}$ by%
\[
S=S(\delta_{0})=\{z_{0}\in V_{H}^{2r}:z_{0}=u_{0}-v_{n}^{\pm},\ \left\Vert
z\left(  t\right)  \right\Vert _{V_{H}^{2r}}\leq\delta_{0}\text{ }\forall
t\geq0\},
\]%
\[
U=U(\delta_{0})=\left\{
\begin{array}
[c]{c}%
z_{0}\in V_{H}^{2r}:z_{0}=u_{0}-v_{n}^{\pm}\text{ and there is a complete
trajectory }\phi\text{ }\\
\text{such that }\phi\left(  0\right)  =u^{0}\text{ and }\left\Vert
\phi\left(  t\right)  -v_{n}^{\pm}\right\Vert _{V_{H}^{2r}}\leq\delta
_{0}\text{ }\forall t\leq0
\end{array}
\right\}  .
\]
We observe that if $z_{0}\in S$, then $u\left(  t\right)
\overset{t\rightarrow+\infty}{\rightarrow}v_{n}^{\pm}$ in $V_{H}^{2r}$. This
follows from the fact that any solution $u\left(  \text{\textperiodcentered
}\right)  $ converges forwards in time to a fixed point in $L^{2}(0,1)$ (see
\cite[Lemma 5.1]{ArrRBVal}) and that using inequality (\ref{IneqEta}) for
$u\left(  \text{\textperiodcentered}\right)  $ we deduce that the convergence
is true in the space $V_{H}^{2r}$ as well. Also, by the properties of the
bounded complete trajectories given in the previous section, we infer that if
$z_{0}\in U$, then $\phi\left(  t\right)  \overset{t\rightarrow-\infty
}{\rightarrow}v_{n}^{\pm}$ in $V_{H}^{2r}.$

Let $X_{i}^{2\alpha}=V_{i}^{\prime}\cap V^{2\alpha}$, where $3/4<\alpha<1 $,
so $V^{2\alpha}=X_{1}^{2\alpha}\oplus X_{2}^{2\alpha}$.

\bigskip

\begin{theorem}
$S\left(  \delta_{0}\right)  $ is tangent to $X_{2}^{2\alpha}$ at the origin
as $\delta_{0}\rightarrow0$. $U\left(  \delta_{0}\right)  $ is tangent to
$X_{1}^{2\alpha}$ at the origin as $\delta_{0}\rightarrow0$.
\end{theorem}

\begin{proof}
Let $z_{0}\in S$. As before, we split the function $z\left(
\text{\textperiodcentered}\right)  $ by $z\left(  t\right)  =z_{1}%
(t)+z_{2}(t)$, where $z_{i}\left(  \text{\textperiodcentered}\right)  $
satisfy (\ref{VC}). Using (\ref{IneqL2D}), $\left\Vert E_{1}v\right\Vert
_{V^{2\alpha}}\leq C_{1}\left\Vert v\right\Vert _{V^{2\alpha}}$ and
$V_{H}^{2r}\subset V^{2\alpha}$ we obtain%
\[
\left\Vert e^{L_{1}t}z_{1}(t)\right\Vert _{V^{2\alpha}}\leq C_{2}%
e^{-\lambda_{m}t}\left\Vert z_{1}(t)\right\Vert _{V^{2\alpha}}\leq
C_{3}e^{-\lambda_{m}t}\left\Vert z(t)\right\Vert _{V_{H}^{2r}}\leq C_{3}%
\rho_{0}e^{-\lambda_{m}t}\overset{t\rightarrow+\infty}{\rightarrow}0.
\]
Hence, from the equality%
\[
e^{L_{1}t}z_{1}(t)=z_{1}\left(  0\right)  +\int_{0}^{t}e^{L_{1}s}E_{1}g\left(
z(s)\right)  ds
\]
we get%
\[
E_{1}z_{0}=z_{1}\left(  0\right)  =-\int_{0}^{\infty}e^{L_{1}s}E_{1}g\left(
z(s)\right)  ds.
\]%
\begin{align*}
z(t)  & =e^{-L_{2}t}z_{2}\left(  0\right)  +\int_{0}^{t}e^{-L_{2}\left(
t-s\right)  }E_{2}g\left(  z(s)\right)  ds-\int_{0}^{\infty}e^{-L_{1}%
(t-s)}E_{1}g\left(  z(s)\right)  ds+\int_{0}^{t}e^{-L_{1}(t-s)}E_{1}g\left(
z(s)\right)  ds\\
& =e^{-L_{2}t}z_{2}\left(  0\right)  +\int_{0}^{t}e^{-L_{2}\left(  t-s\right)
}E_{2}g\left(  z(s)\right)  ds-\int_{t}^{\infty}e^{-L_{1}(t-s)}E_{1}g\left(
z(s)\right)  ds.
\end{align*}
Thus, on the one hand, using (\ref{IneqL2C}) and (\ref{Acotz2Gen}) we have%
\begin{align}
\left\Vert E_{1}z_{0}\right\Vert _{V^{2\alpha}}  & \leq\int_{0}^{\infty
}\left\Vert e^{L_{1}s}E_{1}g\left(  z(s)\right)  \right\Vert _{V^{2\alpha}%
}ds\leq C_{4}\left\Vert E_{1}\right\Vert \int_{0}^{\infty}e^{\lambda_{m}%
s}\left\Vert g\left(  z(s)\right)  \right\Vert _{H^{-1}}ds\nonumber\\
& \leq C_{5}\int_{0}^{\infty}e^{\lambda_{m}s}\left\Vert z(s)\right\Vert
_{V_{H}^{2r}}^{\frac{1}{2}}\left\Vert z(s)\right\Vert _{V^{2\alpha}}ds\leq
C_{6}\rho_{0}^{1/2}\sup_{t\geq0}\left\Vert z(t)\right\Vert _{V^{2\alpha}%
}.\label{E1z0}%
\end{align}
On the other hand, by (\ref{IneqL2A})-(\ref{IneqL2C}) and (\ref{Acotz2Gen}) we
have
\begin{align*}
\left\Vert z(t)\right\Vert _{V^{2\alpha}}  & \leq C_{7}e^{-\lambda_{m+1}%
t}\left\Vert E_{2}z_{0}\right\Vert _{V^{2\alpha}}+C_{8}\left\Vert
E_{2}\right\Vert \int_{0}^{t}e^{-\lambda_{m+1}(t-s)}(t-s)^{-\alpha}\left\Vert
z(s)\right\Vert _{V_{H}^{2r}}^{\frac{1}{2}}\left\Vert z(s)\right\Vert
_{V^{2\alpha}}ds\\
& +C_{9}\left\Vert E_{1}\right\Vert \int_{t}^{\infty}e^{-\lambda_{m}%
(t-s)}\left\Vert z(s)\right\Vert _{V_{H}^{2r}}^{\frac{1}{2}}\left\Vert
z(s)\right\Vert _{V^{2\alpha}}ds,
\end{align*}
so%
\begin{align*}
\sup_{t\geq0}\left\Vert z(t)\right\Vert _{V^{2\alpha}}  & \leq C_{7}\left\Vert
E_{2}z_{0}\right\Vert _{V^{2\alpha}}+\delta_{0}^{1/2}C_{10}\left(  \int%
_{0}^{\infty}e^{-\lambda_{m+1}s}s^{-\alpha}ds+\int_{-\infty}^{0}%
e^{-\lambda_{m}s}ds\right)  \sup_{t\geq0}\left\Vert z(t)\right\Vert
_{V^{2\alpha}}\\
& =C_{7}\left\Vert E_{2}z_{0}\right\Vert _{V^{2\alpha}}+\delta_{0}^{1/2}%
C_{11}\sup_{t\geq0}\left\Vert z(t)\right\Vert .
\end{align*}
If $\delta_{0}^{1/2}C_{11}\leq1/2$, then%
\[
\sup_{t\geq0}\left\Vert z(t)\right\Vert _{V^{2\alpha}}\leq2C_{7}\left\Vert
E_{2}z_{0}\right\Vert _{V^{2\alpha}}.
\]
Therefore, by (\ref{E1z0})\ we obtain%
\[
\left\Vert E_{1}z_{0}\right\Vert _{V^{2\alpha}}\leq C_{12}\rho_{0}%
^{1/2}\left\Vert E_{2}z_{0}\right\Vert _{V^{2\alpha}},
\]
proving that $S\left(  \delta_{0}\right)  $ is tangent to $X_{2}^{2\alpha}$ at
the origin as $\delta_{0}\rightarrow0$.

Let $z_{0}\in U$. Denote $\phi_{i}(t)=E_{i}\phi\left(  t\right)  $ and
$z_{i}(t)=\phi_{i}(t)-v_{n}^{\pm},\ z(t)=z_{1}(t)+z_{2}(t)$ for $t\in
\mathbb{R}$. Using (\ref{IneqL2B}), $\left\Vert E_{2}v\right\Vert
_{V^{2\alpha}}\leq D_{1}\left\Vert v\right\Vert _{V^{2\alpha}}$ and
$V_{H}^{2r}\subset V^{2\alpha}$ we obtain%
\[
\left\Vert e^{-L_{2}t}z_{2}(-t)\right\Vert _{V^{2\alpha}}\leq D_{2}%
e^{-\lambda_{m+1}t}\left\Vert z_{2}(-t)\right\Vert _{V^{2\alpha}}\leq
D_{3}e^{-\lambda_{m+1}t}\left\Vert z_{2}(-t)\right\Vert _{V_{H}^{2r}}\leq
D_{3}\rho_{0}e^{-\lambda_{m+1}t}\overset{t\rightarrow+\infty}{\rightarrow}0.
\]
Hence, (\ref{VC}) gives for $t>0$ that%
\begin{align*}
E_{2}z_{0}  & =z_{2}(0)=e^{-L_{2}t}z_{2}(-t)+\int_{0}^{t}e^{-L_{2}(t-s)}%
E_{2}g(z(s-t))ds\\
& =e^{-L_{2}t}z_{2}(-t)+\int_{-t}^{0}e^{L_{2}s}E_{2}%
g(z(s))ds\overset{t\rightarrow+\infty}{\rightarrow}\int_{-\infty}^{0}%
e^{L_{2}s}E_{2}g(z(s))ds.
\end{align*}
Hence, for $t\leq0$ we get%
\begin{align*}
z(t)  & =e^{-L_{1}t}z_{1}(0)+\int_{0}^{t}e^{-L_{1}(t-s)}E_{1}g(z(s))ds+\int%
_{-\infty}^{0}e^{-L_{2}(t-s)}E_{2}g(z(s))ds-\int_{t}^{0}e^{-L_{2}(t-s)}%
E_{2}g(z(s))ds\\
& =e^{-L_{1}t}z_{1}(0)+\int_{0}^{t}e^{-L_{1}(t-s)}E_{1}g(z(s))ds+\int%
_{-\infty}^{t}e^{-L_{2}(t-s)}E_{2}g(z(s))ds.
\end{align*}
By (\ref{IneqL2A}) and (\ref{Acotz2Gen}) we deduce that%
\begin{align}
\left\Vert E_{2}z_{0}\right\Vert _{V^{2\alpha}}  & \leq\int_{-\infty}%
^{0}\left\Vert e^{L_{2}s}E_{2}g\left(  z(s)\right)  \right\Vert _{V^{2\alpha}%
}ds\leq D_{4}\left\Vert E_{2}\right\Vert \int_{-\infty}^{0}e^{\lambda_{m+1}%
s}\left(  -s\right)  ^{-\alpha}\left\Vert g\left(  z(s)\right)  \right\Vert
_{H^{-1}}ds\nonumber\\
& \leq D_{5}\int_{-\infty}^{0}e^{\lambda_{m+1}s}\left(  -s\right)  ^{-\alpha
}\left\Vert z(s)\right\Vert _{V_{H}^{2r}}^{\frac{1}{2}}\left\Vert
z(s)\right\Vert _{V^{2\alpha}}ds\leq D_{6}\rho_{0}^{1/2}\sup_{t\leq
0}\left\Vert z(t)\right\Vert _{V^{2\alpha}}.\label{E2z0}%
\end{align}
Also, (\ref{IneqL2A}), (\ref{IneqL2C})-(\ref{IneqL2D}) and (\ref{Acotz2Gen})
imply for $t\leq0$ that%
\begin{align*}
\left\Vert z(t)\right\Vert _{V^{2\alpha}}  & \leq D_{7}e^{-\lambda_{m}%
t}\left\Vert E_{1}z_{0}\right\Vert _{V^{2\alpha}}+D_{8}\left\Vert
E_{1}\right\Vert \int_{t}^{0}e^{-\lambda_{m}(t-s)}\left\Vert z(s)\right\Vert
_{V_{H}^{2r}}^{\frac{1}{2}}\left\Vert z(s)\right\Vert _{V^{2\alpha}}ds\\
& +D_{9}\left\Vert E_{2}\right\Vert \int_{-\infty}^{t}e^{-\lambda_{m+1}%
(t-s)}(t-s)^{-\alpha}\left\Vert z(s)\right\Vert _{V_{H}^{2r}}^{\frac{1}{2}%
}\left\Vert z(s)\right\Vert _{V^{2\alpha}}ds,
\end{align*}
so%
\begin{align*}
\sup_{t\leq0}\left\Vert z(t)\right\Vert _{V^{2\alpha}}  & \leq D_{7}%
e^{-\lambda_{m}t}\left\Vert E_{1}z_{0}\right\Vert _{V^{2\alpha}}+\delta
_{0}^{1/2}D_{10}\left(  \int_{-\infty}^{0}e^{-\lambda_{m}s}ds+\int_{0}%
^{\infty}e^{-\lambda_{m+1}s}s^{-\alpha}ds\right)  \sup_{t\leq0}\left\Vert
z(t)\right\Vert _{V^{2\alpha}}\\
& =D_{7}\left\Vert E_{1}z_{0}\right\Vert _{V^{2\alpha}}+\delta_{0}^{1/2}%
D_{11}\sup_{t\leq0}\left\Vert z(t)\right\Vert .
\end{align*}
If $\delta_{0}^{1/2}D_{11}\leq1/2$, then%
\[
\sup_{t\leq0}\left\Vert z(t)\right\Vert _{V^{2\alpha}}\leq2D_{7}\left\Vert
E_{1}z_{0}\right\Vert _{V^{2\alpha}}.
\]
Hence, (\ref{E2z0}) implies%
\[
\left\Vert E_{2}z_{0}\right\Vert _{V^{2\alpha}}\leq D_{12}\rho_{0}%
^{1/2}\left\Vert E_{1}z_{0}\right\Vert _{V^{2\alpha}},
\]
so $U\left(  \delta_{0}\right)  $ is tangent to $X_{1}^{2\alpha}$ at the
origin as $\delta_{0}\rightarrow0$.
\end{proof}

\bigskip

We observe that in the neighborhood $O_{\delta_{0}}(v_{n}^{\pm})$ there are
three types of points:

\begin{enumerate}
\item If $z_{0}=u_{0}-v_{n}^{\pm}\in S\backslash\{0\}$, then $u\left(
t\right)  \in O_{\delta_{0}}(v_{n}^{+})$, for $t\geq0$, and $u\left(
t\right)  \overset{t\rightarrow+\infty}{\rightarrow}v_{n}^{\pm}$, being
$z\left(  t\right)  =u(t)-$ tangent to $X_{2}^{2\alpha}$ as $t\rightarrow
+\infty$. Also, if there is a complete trajectory $\phi$ passing through
$u_{0}$, then there is $t_{1}<0$ such that $\phi\left(  t_{1}\right)
\not \in O_{\delta_{0}}(v_{n}^{\pm}).$

\item If $z_{0}=u_{0}-v_{n}^{\pm}\in U\backslash\{0\}$, then there is a
complete trajectory $\phi$ passing through $u_{0}$ such that $\phi\left(
t\right)  \overset{t\rightarrow-\infty}{\rightarrow}v_{n}^{\pm}$, being
$z\left(  t\right)  $ tangent to $X_{1}^{2\alpha}$ as $t\rightarrow-\infty$,
and there is $t_{1}>0$ such that $\phi\left(  t_{1}\right)  \not \in
O_{\delta_{0}}(v_{n}^{\pm}).$

\item If $u_{0}\in O_{\delta_{0}}(v_{n}^{\pm}),\ z_{0}\not \in \{S\cup U\}$,
then there is $t_{1}>0$ such that $u\left(  t_{1}\right)  \not \in
O_{\delta_{0}}(v_{n}^{\pm})$. Moreover, if there is a complete trajectory
$\phi$ passing through $u_{0}$, then there is $t_{2}<0$ such that $\phi\left(
t_{2}\right)  \not \in O_{\delta_{0}}(v_{n}^{\pm}).$
\end{enumerate}

All in all, this means that the fixed point $v_{n}^{\pm}$ satisfies the
saddle-point property. Hence, the fixed point is hyperbolic.

\end{document}